# Preliminary design of Debris removal missions by means of simplified models for Low-Thrust, many-revolution transfers

Federico Zuiani
*PhD Candidate, Space Advance Research Team,*
*School of Engineering, University of Glasgow, James Watt South Building, Glasgow, G12 8QQ, United Kingdom*
Tel.: +44-141-5484558
f.zuiani.1@research.gla.ac.uk

Massimiliano Vasile
*Reader, Space Advance Research Team,*
*Department of Mechanical & Aerospace Engineering, University of Strathclyde, 75 Montrose Street, Glasgow, G1 1XJ, United Kingdom*
Tel.: +44-141-3306465
Fax: +44-141-3305560
massimiliano.vasile@strath.ac.uk

***Abstract***

This paper presents a novel approach for the preliminary design of Low-Thrust (LT), many-revolution transfers. The main feature of the novel approach is a considerable reduction in the control parameters and a consequent gain in computational speed. Each spiral is built by using a predefined pattern for thrust direction and switching structure. The pattern is then optimised to minimise propellant consumption and transfer time. The variation of the orbital elements due to the propulsive thrust is computed analytically from a first-order solution of the perturbed Keplerian motion. The proposed approach allows for a realistic estimation of $\Delta V$ cost and time of flight required to transfer a spacecraft between two arbitrary orbits. Eccentricity and plane changes are both taken into account. The novel approach is applied here to the design of missions for the removal of space debris by means of an Ion Beam Shepherd (IBS) Spacecraft. In particular, two slightly different variants of the proposed low-thrust control model are used for the two main phases of the debris removal mission, i.e. the rendezvous with the target object and its removal. Thanks to their relatively low computational cost they can be included in a multiobjective optimisation problem in which the sequence and timing of the removal of five hypothetical pieces of debris are optimised in order to minimise both propellant consumption and mission duration.

## 1. Introduction

One of the most critical issues related to the exploitation of Space around the Earth is the threat posed by space debris. Since the beginning of the space era in the late 1950s, an increasing number of man-made, inert objects has been orbiting the Earth. Recent statistics revealed around 15000 trackable objects, for a total of some 6000 tons of material. Some of these objects are simply spent upper stages of launch vehicles, some others are satellites which are no longer active due to failures or to having reached their end of life. Others, however, are the results of past collisions. It is easy to imagine that even a single collision between two objects is likely to generate tens of smaller objects as a result. The outcome of a collision in an already crowded environment could generate a cascade of collisions generating an exponentially increasing volume of space debris. In fact, the debris produced by a collision is itself likely to collide with other objects, thereby producing other debris which will generate further collisions, and so on. This chain reaction, known as the *Kessler Syndrome* [1], occurs once the rate of generation of debris due to collisions or simple human-driven additions, exceeds the natural debris removal rate. According to Kessler, this reaction is likely to be ignited once the object density in a certain orbital band reaches a critical point; once started, it will probably render most spacecraft in that orbital band useless within a matter of months or years.

Recent guidelines issued by international spacer regulatory institutions such as the United Nations Committee for the Peaceful Uses of Outer Space (COPUOS) [2] and the Inter-Agency Space Debris Coordination Committee (IADC) [3] prescribe some actions to be followed by national or private agencies putting satellites into orbit in order to mitigate debris growth. For example, it is demanded that every new mission in Low Earth Orbit (LEO) must be planned such that the satellite itself must re-enter in the Earth's atmosphere within 25 years after the end of the mission. Alternatively, for higher orbits like Geostationary orbits, the requirement is for the spacecraft to be placed on a higher *graveyard* orbit. Measures like these, even if strictly applied (and at the moment compliance with them is on a voluntary basis) are only likely to slow down the accumulation of space debris around the Earth. Therefore, active removal actions will probably be needed in the near future to eliminate at least the most dangerous objects.

There have been various proposals on how to remove inert objects from space. They can be generally classified in two major groups: contactless and with direct physical contact. In the latter category one can find methods based on some form of docking with or capturing the object. Once the removing spacecraft and the piece of debris are attached, the latter is dragged into a re-entry trajectory or to a graveyard orbit. Technical problems related to the attitude state of motion of the piece of debris and the fragility of appendices and cover material (including paint) make this removal solution complicated. A potentially interesting solution is represented by Project ROGER [4], developed by EADS/Astrium with the support of ESA. Among contactless solutions on can find what is commonly referred to as the *space broom* [5]. It entails irradiating the target object with a high-power laser which will induce sublimation of the surface material; the ejecta plume will then generate a low thrusting acceleration which will slowly degrade the debris' orbit until it reaches an altitude where atmospheric drag will accelerate its re-entry. Such a technique has the advantage that no physical contact is required, on the other hand current proposals envisage the use of lasers installed on Earth and beaming through the atmosphere. The beam collimation and thrust time is therefore limited and this solution is effective for small-sized objects only. Recent proposals have demonstrated that the use of in-space lasers systems might be more interesting even to remove larger objects [6]. Other proposals involve for example the use of electrodynamic tethers [7], inflatable balloons [8], which are meant to be lightweight and efficient but require, however, the physical attachment of the device to the target object and are therefore of difficult application to existing debris.

A recent idea simultaneously proposed by Bombardelli et al. [9], Bonnal et al. [10] and JAXA [**Error! Reference source not found.**] suggested the use of a collimated beam of ions generated by a spacecraft flying in formation with the piece of debris. In this paper this concept will be called Ion Beam Shepherd (IBS), using the name introduced by Bombardelli et al.. The effect of the ion beam is that of producing a thrusting force, equal in magnitude but opposite in direction, on both the IBS and the piece of debris. This force will induce a thrusting acceleration which can be controlled in order to modify the orbit of the piece of debris. A second ion engine is then fired in a direction opposite to the first one in order to keep the IBS spacecraft at a constant distance from the piece of debris. Among the advantages of this concept is the fact that it employs already existing and proven technologies; it does not require any contact with the target, and the fact that a single spacecraft can be used to fetch and deorbit multiple pieces of debris. In [6] one can find a similar concept that uses concentrated solar light instead of ions to generate a thrust and modify the orbit of debris.

Assuming a scenario in which a single IBS needs to de-orbit multiple pieces of debris, one would need to solve an interesting mission design problem: the optimisation of the de-orbit sequence and trajectories for multiple target objects in minimum time and with minimum propellant. In the hypothetical mission scenario which is analysed in this work, it is assumed that a number of pieces of debris have been shortlisted as priority targets due to the threat they pose to satellites operating in LEO. For example Johnson et al. [11] propose some criteria to choose the object whose removal will be most effective to mitigate the risk of collisions. They underline that an effective removal strategy must be targeted first to large objects in crowded orbits up to 1500 km. Thus , a removal mission by means of an IBS spacecraft is planned to be launched from the Earth. Its task is that of removing five objects lying on different low Earth orbits. The design of such a mission is a complex optimisation problem, because it requires the computation of multiple low-thrust, many-revolution transfers. Therefore, this work proposes an approach to the fast estimation and optimisation of the cost and time duration of the fetch and de-orbit sequences. In past works, other authors have already proposed approaches to the design of low thrust, many revolution transfers, based on analytical solutions to an optimal orbit raising problem under the assumption of small eccentricity [12,13,14] or on averaging techniques [15,16,17]. This work, proposes a different approach, based on a first order solution of perturbed Keplerian motion. The approach in this paper aims at capturing the definition and optimisation of the thrusted arcs for each orbit without sacrificing computational speed. The approach can be classified as direct method for trajectory design as it does not derive the necessary conditions for optimality but translates the initial optimal control problem into an NLP problem.

The paper is organised as follows: Section 2 will briefly outline the IBS concept and in particular will outline how to compute the thrusting acceleration generated on a given target object; Section 3 will analyse an hypothetical mission profile for the removal mission and most important, Subsections 3.1 and 3.2 will present in detail the proposed trajectory models. Section 4 will then show how the mission design problem can be then translated into as a series of multiobjective optimisation problems which are solved with a stochastic optimiser. The results are then presented and discussed.

## 2. The Ion Beam Shepherd Spacecraft

As shown by Bombardelli et al. [9], the concept behind the Ion Beam Shepherd is relatively simple and envisions employing a spacecraft provided with two sets of Ion engines mounted along the same axis but in opposite directions (see Fig. 1). The jet from one of the sets will be directed towards the piece of debris and will exert a thrusting force $\mathbf{F}_{p1}$ on it. Due to Newton's third law, an opposite force of same magnitude will also act on the spacecraft itself, but this component will be balanced by the thrust $\mathbf{F}_{p2}$ provided by the other set of Ion engines.

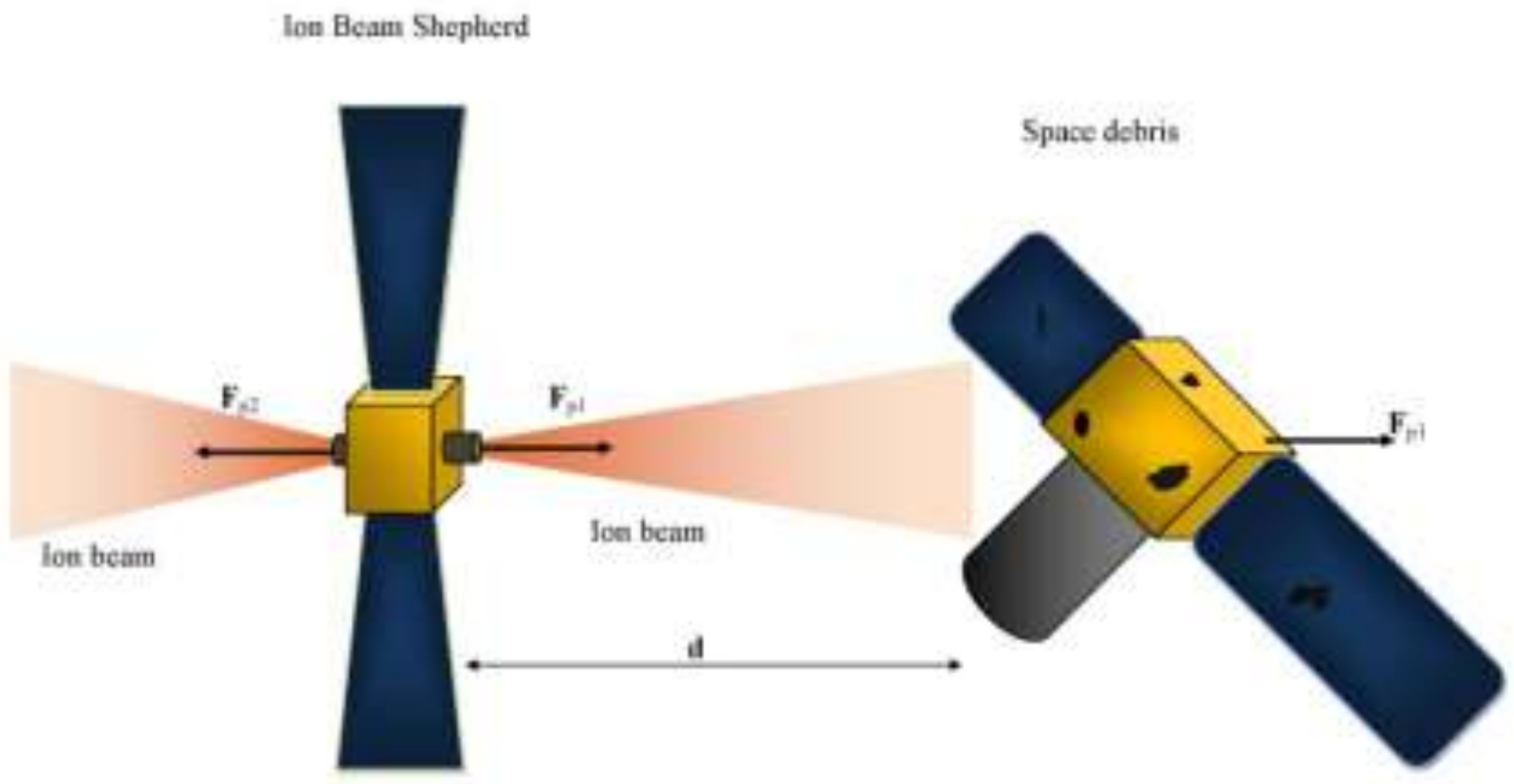

**Fig. 1** Ion Beam Shepherd spacecraft.

Since it is necessary to keep the Shepherd spacecraft at a constant distance from the debris, the thrust $\mathbf{F}_{p2}$ shall be such that the second derivative of the distance d between the two spacecraft is null:

$$\ddot{\mathbf{d}} = \frac{\mathbf{F}_{p2} - \mathbf{F}_{p1}}{m_{IBS}} - \frac{\mathbf{F}_{p1}}{m_d} = 0 \tag{1}$$

Note that in Eq. (1) the acceleration terms due to the gravity of the central body have been neglected since it is assumed that the debris and the Shepherd are in close proximity and arranged in a leader-follower configuration. A more accurate and detailed modelling of the proximal motion dynamics of these two bodies is beyond the scope of this study. Thus in the following sections, the IBS-debris combination will be treated as a point mass, in order to apply two-body dynamics. By rearranging the terms in Eq. (1) one obtains:

$$F_{p2} = F_{p1}\left(1 + \frac{m_{IBS}}{m_d}\right) \tag{2}$$

Under the assumption that the total propulsive power of the IBS spacecraft $P_{tot}$ is constant and that the total propulsive thrust is proportional to it $F_{tot}$, one can write:

$$F_{p1} + F_{p2} = F_{tot} \propto P_{tot} \tag{3}$$

thus:

$$F_{tot} = F_{p1}\left(2 + \frac{m_{IBS}}{m_d}\right) \tag{4}$$

Therefore, the maximum acceleration acting on the IBS-debris combination can be computed as a function of the total available thrust $F_{tot}$:

$$\varepsilon_{IBS-debr} = \frac{F_{p1}}{m_d} = \frac{F_{tot}}{2m_d + m_{IBS}} \tag{5}$$

It is assumed here to have a 1000 kg IBS spacecraft with a total available thrust of 0.5 N. Such a high thrust would correspond to a substantial power and propulsion system mass, however this is deemed realistic if one considers that the payload of the IBS spacecraft is in fact its propulsion and power systems. Hence, the propulsion and power systems might be oversized compared to other applications in which ion engines are used for propulsion only. Note that the validity of the methodology proposed in this paper would not be affected even if lower thrust levels were considered. Thus, in this case, considering for example an 800 kg debris, the magnitude of the acceleration, would be $1.923 \cdot 10^{-7}$ km/s$^2$. If one considers instead the spacecraft alone, the acceleration achievable would be slightly higher, $5 \cdot 10^{-7}$ km/s$^2$. Given this order of magnitude, the thrust acceleration can be considered as a perturbative force compared to the Earth's gravitational force and therefore the analytical approach to the propagation of the low-thrust motion described in [18] can be applied.

## 3. Mission profile

The objective of this study is that of optimising the performance and cost of a debris de-orbiting mission performed by a single spacecraft. As mentioned in the introduction, it is assumed that there are five pieces of debris of different masses and lying in circular orbits with different radii and orientations. It is assumed that, the IBS spacecraft departs from a low-Earth parking orbit, rendezvous with the first object, transfer it to an elliptical re-entry orbit, rendezvous with the second object, transfers it to a second elliptical re-entry orbit, and so on and so forth until all five pieces of debris are removed. One important issue is defining in which order the pieces of debris need to be de-orbited. In the following all possible sequences are generated a priori and optimised one by one.

Each fetch and de-orbit operation is split in two phases:

- A *de-orbit* phase, in which the perigee of the orbit of the piece of debris is lowered such that the orbit will decay naturally in a relatively short time. In this study it is assumed that this condition is met if the perigee altitude of the debris' orbit is equal or lower than 300 km.
- A *transfer* phase, in which the IBS spacecraft rendezvouses with the next piece of debris (which lies on a circular orbit), after having abandoned the current piece of debris on an orbit with a 300 km perigee altitude.

Given the magnitude of the available thrust acceleration, both phases require a spiral orbit transfer. If a direct transcription approach is used to optimise each spiral the number of parameters that needs to be defined is very high leading to high computational times. The latter fact would make the solution of a multiobjective optimisation of all possible de-orbiting sequences computationally intractable. Thus, in this paper a simplified, highly efficient, trajectory model is proposed for each one of the two phases.

### 3.1. De-Orbiting Trajectory Model

The objective of the de-orbiting phase is that of lowering an initial circular orbit such that its perigee is equal or below 300 km, which basically translates into a perigee lowering manoeuvre. Therefore, it is appropriate to assume that in general, as soon as the initial circular orbit becomes slightly eccentric, one keeps thrusting around the apogee in order to lower the perigee. The thrust level will also be kept at its maximum in order to minimize gravity losses. Moreover, since the de-orbit condition is independent of the final orbit's orientation, one can reasonably assume that the perigee lowering will be performed in-plane. In this sense, the only Keplerian parameters which need to be altered are the semi-major axis and eccentricity. By analysing the structure of Gauss' variational equations:

$$
\begin{aligned}
\frac{da}{dt} &= \frac{2a^2}{h}\left(e\sin\theta a_r + \frac{p}{r}a_\theta\right) \\
\frac{de}{dt} &= \frac{1}{h}\left\{p\sin\theta a_r + \left[\left(p+r\right)\cos\theta + re\right]a_\theta\right\} \\
\frac{di}{dt} &= \frac{r\cos\theta}{h}a_h \\
\frac{d\Omega}{dt} &= \frac{r\sin\theta}{h\sin i}a_h \\
\frac{d\omega}{dt} &= \frac{1}{he}\left[-p\cos\theta a_r + \left(p+r\right)\sin\theta a_\theta\right] - \frac{r\sin\theta\cos i}{h\sin i}a_h \\
\frac{d\theta}{dt} &= \frac{h}{r^2} + \frac{1}{eh}\left[p\cos\theta a_r - \left(p+r\right)\sin\theta a_\theta\right]
\end{aligned}
\tag{6}
$$

where $a$ is the semi-major axis, $e$ the eccentricity, $i$ the inclination $\Omega$ is the Right Ascension of the Ascending node (RAAN), $\omega$ the argument of perigee, $\theta$ the true anomaly, $p$ the semi-latus rectus and $h$ the angular momentum; $a_r$, $a_\theta$ and $a_h$ are the radial, transversal and out-of-plane components of the thrust acceleration. If one considers the case of thrusting with maximum acceleration along arcs which are symmetrical around apogee ($\theta=\pi$), one can see that the contributions to semi-major axis and eccentricity variations given by the $a_r$ components are negligible (since they are multiplied by $\sin\theta$). Therefore, a good suboptimal thrust direction can be obtained by imposing $a_\theta$ as the only non-zero component of the thrust acceleration. Under these assumptions the variation of Keplerian parameters will be [19]:

$$
\begin{aligned}
\frac{da}{dt} &= \frac{2a^2}{h}\frac{p}{r}a_\theta \\
\frac{de}{dt} &= \frac{1}{h}\left[\left(p+r\right)\cos\theta + re\right]a_\theta \\
\frac{di}{dt} &= 0 \\
\frac{d\Omega}{dt} &= 0 \\
\frac{d\omega}{dt} &= \frac{\left(p+r\right)}{he}\sin\theta a_\theta \\
\frac{d\theta}{dt} &= \frac{h}{r^2} + \left(p+r\right)\sin\theta a_\theta
\end{aligned}
\tag{7}
$$

It should be noted that the terms of the variation of $\omega$ and $\theta$ which depend on $a_\theta$ will also be very small due to the presence of $\sin\theta$ integrated around $\theta=\pi$. Fig. 2 visualises the proposed pattern for thrusting arcs.

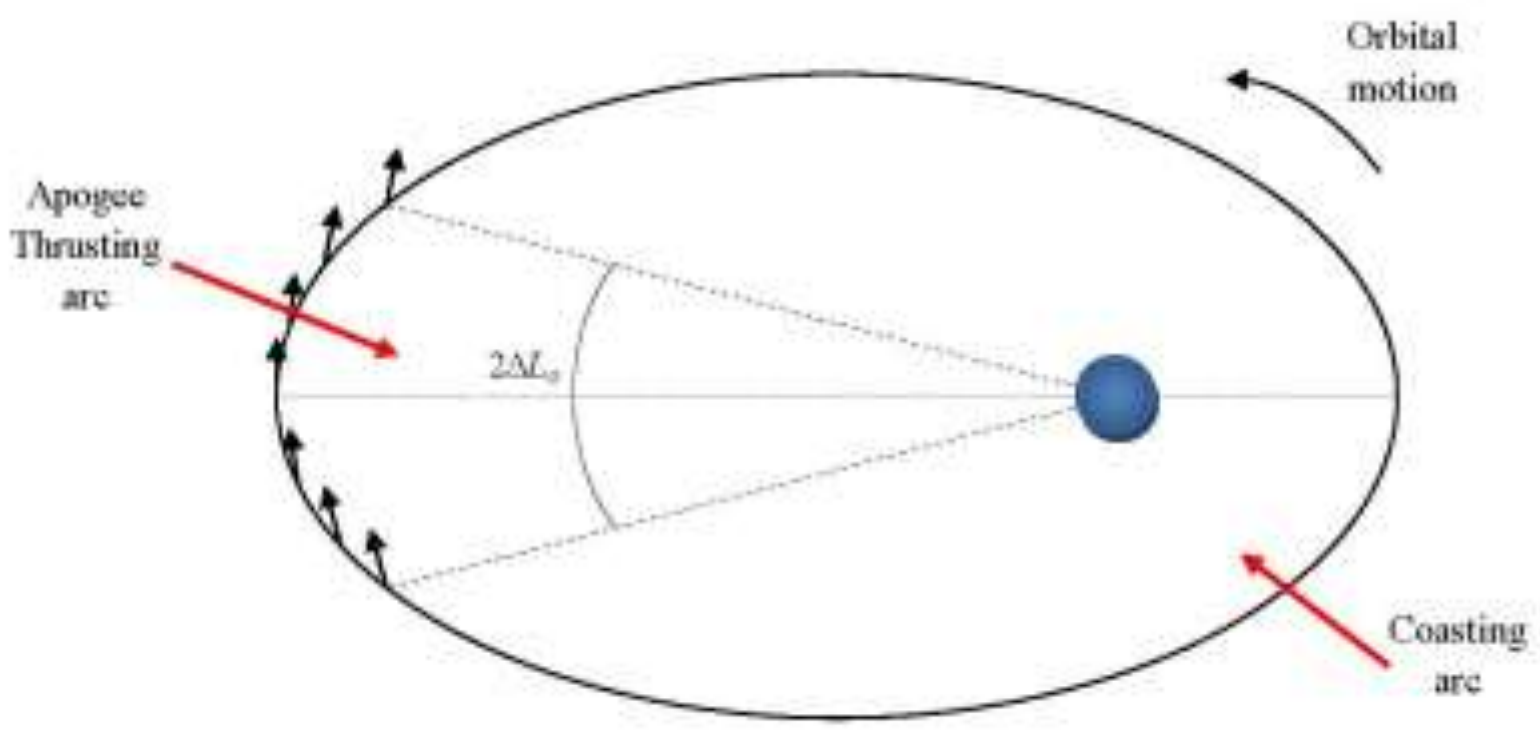

**Fig. 2** Thrusting arc around apogee with thrust directed along transverse direction

In order to obtain a fast propagation of the thrusting arcs, the analytical propagation of perturbed motion with Finite Perturbative Elements in Time (FPET) derived in [18,20] will be used. In order to employ FPET, one has also to assume that the thrust acceleration is constant around each thrusting arc, which is reasonable given the low propellant consumption per arc. This assumption ensures that, if the engine thrust is constant, the resulting acceleration can also be considered constant over short thrusting arcs.

Motion propagation with FPET is based on a first-order analytical solution of perturbed Keplerian motion. In this formulation, the state is expressed in non-singular equinoctial elements:

$$\begin{gathered} a \\ P_1 = e\sin\left(\Omega+\omega\right) \\ P_2 = e\cos\left(\Omega+\omega\right) \\ Q_1 = \tan\frac{i}{2}\sin\Omega \\ Q_2 = \tan\frac{i}{2}\cos\Omega \\ L = \Omega+\omega+\theta \end{gathered} \tag{8}$$

Assuming constant thrust-acceleration in the radial-transverse reference frame:

$$\mathbf{a} = \begin{bmatrix} a_r & a_\theta & a_h \end{bmatrix}^T = \varepsilon\begin{bmatrix} \cos\alpha\cos\beta & \sin\alpha\cos\beta & \sin\beta \end{bmatrix}^T \tag{9}$$

then one can obtain a first-order analytical expansion of the variation of Equinoctial elements, parameterised in Longitude $L$ and with respect to a reference longitude $L_0$:

$$\begin{aligned} a(L) &= a_0(L_0) + \varepsilon\, a_1(L_0, \Delta L, \alpha, \beta) \\ P_1(L) &= P_{10}(L_0) + \varepsilon\, P_{11}(L_0, \Delta L, \alpha, \beta) \\ P_2(L) &= P_{20}(L_0) + \varepsilon\, P_{21}(L_0, \Delta L, \alpha, \beta) \\ Q_1(L) &= Q_{10}(L_0) + \varepsilon\, Q_{11}(L_0, \Delta L, \alpha, \beta) \\ Q_2(L) &= Q_{20}(L_0) + \varepsilon\, Q_{21}(L_0, \Delta L, \alpha, \beta) \\ t(L) &= t_{00}(L_0, \Delta L) + \varepsilon\, t_{11}(L_0, \Delta L, \alpha, \beta) \end{aligned} \tag{10}$$

where:

$$L = L_0 + \Delta L \tag{11}$$

where a $a_0$ $P_{10}$, $P_{20}$, $Q_{10}$, $Q_{20}$ are reference values at $L_0$ and $a_1$, $P_{11}$, $P_{21}$, $Q_{11}$, $Q_{21}$ are first-order terms as reported in [18,20]. In [18] it has also been shown that this analytical propagation scheme provides good accuracy along relatively long trajectory arcs.

As explained above, the only non-zero component of the acceleration will be $a_\theta$ and since the aim is obtaining a decrease of the orbit energy it will also be in the negative direction. Therefore the acceleration azimuth will be $\alpha=-\pi/2$ and the elevation $\beta=0$ (since, as already mentioned, the motion will be within the initial orbit plane). The variation of equinoctial elements after an apogee thrusting arc will be given by:

$$\mathbf{E}^+ = \begin{bmatrix} a \\ P_1 \\ P_2 \end{bmatrix}_{L^+=L_a+\Delta L_a} = \begin{bmatrix} a \\ P_1 \\ P_2 \end{bmatrix}_{L^-=L_a-\Delta L_a} + \varepsilon \begin{bmatrix} a_1\left(2\Delta L_a, \frac{\pi}{2}, 0\right) \\ P_{11}\left(2\Delta L_a, \frac{\pi}{2}, 0\right) \\ P_{21}\left(2\Delta L_a, \frac{\pi}{2}, 0\right) \end{bmatrix} = \mathbf{E}^- + \varepsilon \mathbf{f}\left(2\Delta L_a, \frac{\pi}{2}, 0\right) \tag{12}$$

where $L_a$ is the apocentre longitude, $L^-$ and $L^+$ are the longitudes at the start and end of thrusting respectively. $\Delta L_a$ is the semi-amplitude of the apogee thrusting arc. Note that, given that $\beta=0$, there is no variation on $Q_1$ and $Q_2$ and thus they are omitted. The coasting time is computed from the last of Eqs. (10) as:

$$t_{thrust} = t_{00}\left(L_a - \Delta L_a, 2\Delta L_a\right) + \varepsilon t_{11}\left(2\Delta L_a, -\frac{\pi}{2}, 0\right) \tag{13}$$

Since the thrust magnitude and direction are fixed, the only free control parameter is the semi-amplitude $\Delta L_a$ for each orbit. In order to keep the number of decision variables to a minimum, the semi-amplitude for each orbit is computed from a piece-wise linear polynomial interpolating a limited number of $\Delta L_{a,i}$ over a number of orbits. The nodes $\Delta L_{a,i}$ are equally distributed between orbit 1 and an arbitrary number of orbits (in this paper 1200 was found to be adequate). In this paper the number of interpolating nodes was limited to 2: $\Delta L_{a1}$ and $\Delta L_{af}$.

In order to evaluate the time and $\Delta V$ needed to de-orbit a piece of debris from its initial orbit with semi-major axis $a_{\text{debr0}}$, given a set of decision (or control) parameters $\Delta L_{a1}$ and $\Delta L_{af}$, the following procedure was implemented:

1. Compute the set of initial Equinoctial parameters $L_0$ and $\mathbf{E}_0 = \left[a_{debr0} \quad P_{10} \quad P_{20}\right]^T$ where $P_{10}$ and $P_{20}$ will be null due to the fact that the initial orbit is circular.
2. Initialise the number of orbits, the total $\Delta V$ and time of flight to zero:

$$N_{orbit} = 0$$
$$\Delta V = 0$$
$$ToF = 0$$

3. Set $\mathbf{E}^- = \mathbf{E}_0$ and $L_{coast} = L_0$.
4. Initialise the mass of the IBS spacecraft:

$$m_{IBS} = m_{IBS0}$$

5. While $N_{orbit}$ is smaller than $N_{orbitsmax}$:
   a. $$N_{orbit} = N_{orbit} + 1$$
   b. Interpolate the amplitude of the thrusting arc in the current orbit, i.e. $\Delta L_a$ and compute $L^- = L_a - \Delta L_a$ and $L^+ = L_a + \Delta L_a$.
   c. Compute the acceleration $\varepsilon_{IBS\text{-}debr}$ acting on the IBS-debris combination from Eq. (5).
   d. Compute the time of flight $t_{coast}$ spent coasting from $L_{coast}$ to $L^-$.
   e. Compute the Equinoctial parameters after the thrusting arc $\mathbf{E}^+$ as in Eq. (12).
   f. Compute the current perigee radius $r_p$ and if this is lower than the threshold $\bar{r}_p = 300km$ proceed to step 6, otherwise proceed to step g.

g. Compute the thrusting time $t_{thrust}$ from Eq. (13) and update the total ΔV cost:

$$\Delta V = \Delta V + \varepsilon_{IBS-debr} t_{thrust} \tag{14}$$

h. Update the total time of flight:

$$ToF = ToF + t_{coast} + t_{thrust} \tag{15}$$

i. Update the IBS spacecraft mass:

$$m_{IBS} = \left(m_{IBS} + 2m_{debr}\right)\exp\left(-\frac{\varepsilon_{IBS-debr} t_{thrust}}{I_{sp} g_0}\right) - 2m_{debr} \tag{16}$$

j. Set $\mathbf{E}^{-} = \mathbf{E}^{+}$ and $L_{coast} = L^{+}$

6. Back-track the value of the longitude $L_f$ for which $r_p = \overline{r}_p$ and compute the related and $t_{thrust}$ from Eq. (13) and update $ToF$ and $\Delta V$ accordingly. Compute the Equinoctial parameters $\mathbf{E}_f$ at $L_f$ from Eq. (12).

At this point one gets the $\Delta V$ , the time of flight $ToF$ and the semi-major axis and eccentricity of the final orbit (which are easily computed from $\mathbf{E}_f$). It is important to note that, given the simplifications introduced, once one sets the initial mass and orbit of the piece of debris, and the characteristics of the IBS propulsion system, i.e. $F_{tot}$ and $I_{sp}$, the de-orbit depends exclusively on the mass of the IBS $m_{IBS0}$ at the beginning of the de-orbit phase and the interpolating values for $\Delta L_a$, i.e. $\Delta L_{a1}$ and $\Delta L_{af}$. Therefore, it was decided to pre-compute the corresponding $\Delta V$ and $ToF$ for a given set of these three parameters and for each piece of debris (i.e. for each $m_{debr}$ and $a_{debr0}$). Table 1 reports upper and lower bounds for $m_{IBS0}$, $\Delta L_{a1}$ and $\Delta L_{af}$. and the number of samples taken, equally distributed.

**Table 1** Bounds and number of samples for the de-orbit parameters

| | $m_{IBS0}$ | $\Delta L_{a1}$ | $\Delta L_{af}$ |
|---|---|---|---|
| Lower bound | $m_{dry}$+100=350kg | 0 | 0 |
| Upper bound | $m_{launch}$=1000kg | π | π |
| Samples | 8 | 50 | 50 |

Given the limited number of decision variables, for each piece of debris, one has 20000 de-orbit instances to propagate. Since each instance requires typically $1\cdot10^{-2}$ s of CPU time, with a code implemented in MatLab® and running on a 3.16 GHz, 4 GB desktop PC running Windows 7®, the whole computation can be completed in roughly five minutes. The set of de-orbit $\Delta V$ and $ToF$ is then used to build a response surface, or surrogate model, of the de-orbiting process. Fig. 3a and Fig. 3b show examples of two-dimensional surface, respectively for $\Delta V$ and $ToF$, with respect to a fixed $m_{IBS0}$ of 300 kg. One can see that the two quantities show opposite trends, the $\Delta V$ being high when the $ToF$ is low and vice versa. Fig. 4a and Fig. 4b show the final semi-major axis and eccentricity respectively. Note that the minimum $ToF$ transfer corresponds to a quasi-circular spiralling trajectory in which the IBS spacecraft is thrusting continuously. On the other hand, the minimum $\Delta V$ transfer corresponds also to the one with maximum final eccentricity.

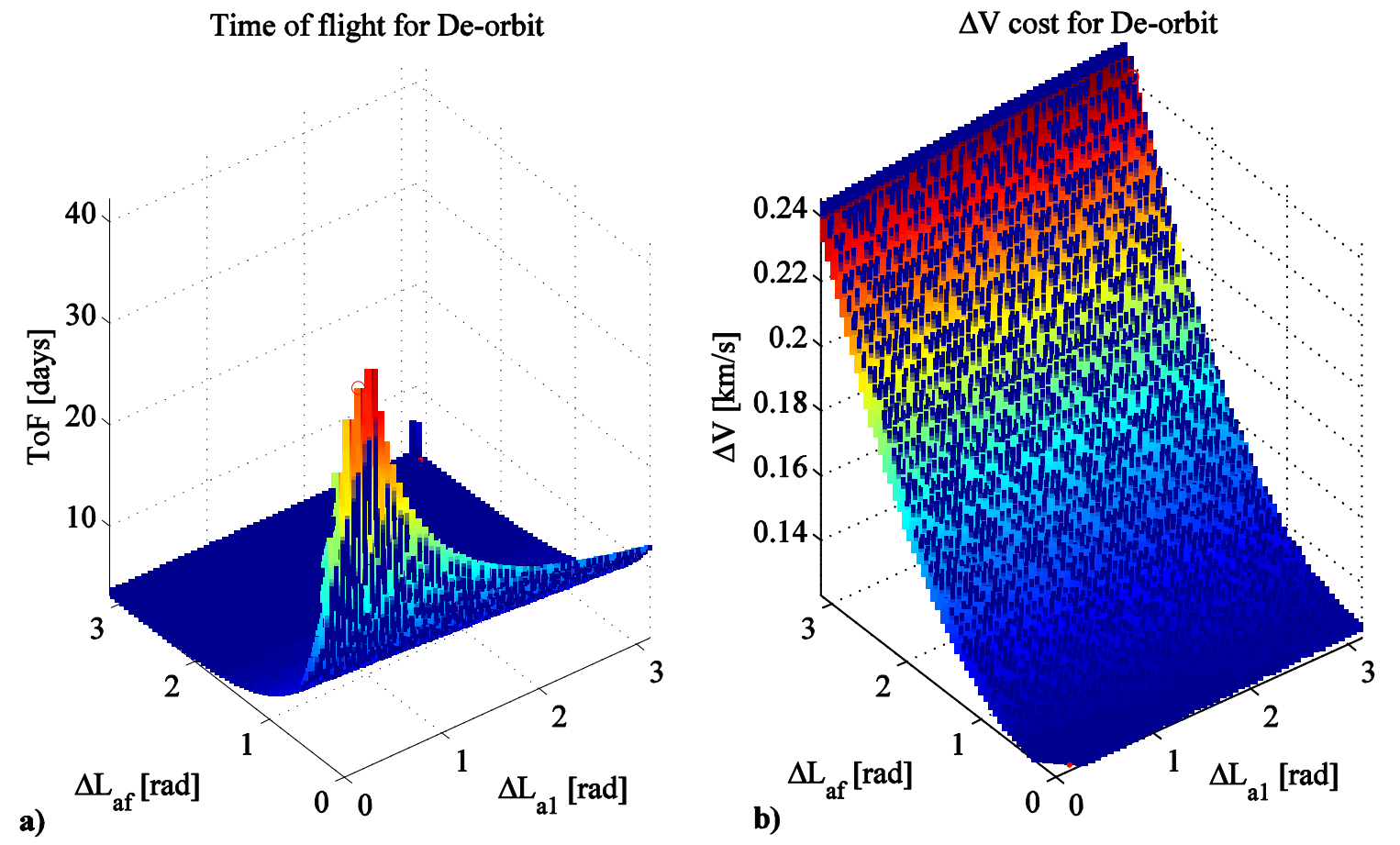

**Fig. 3 a)** $\Delta V$ and **b)** *ToF* surfaces with respect to $\Delta L_{a1}$ and $\Delta L_{af}$ for $m_{IBS0} = 300kg$, $a_{debr0} = 7128km$ and $m_{debr} = 120kg$

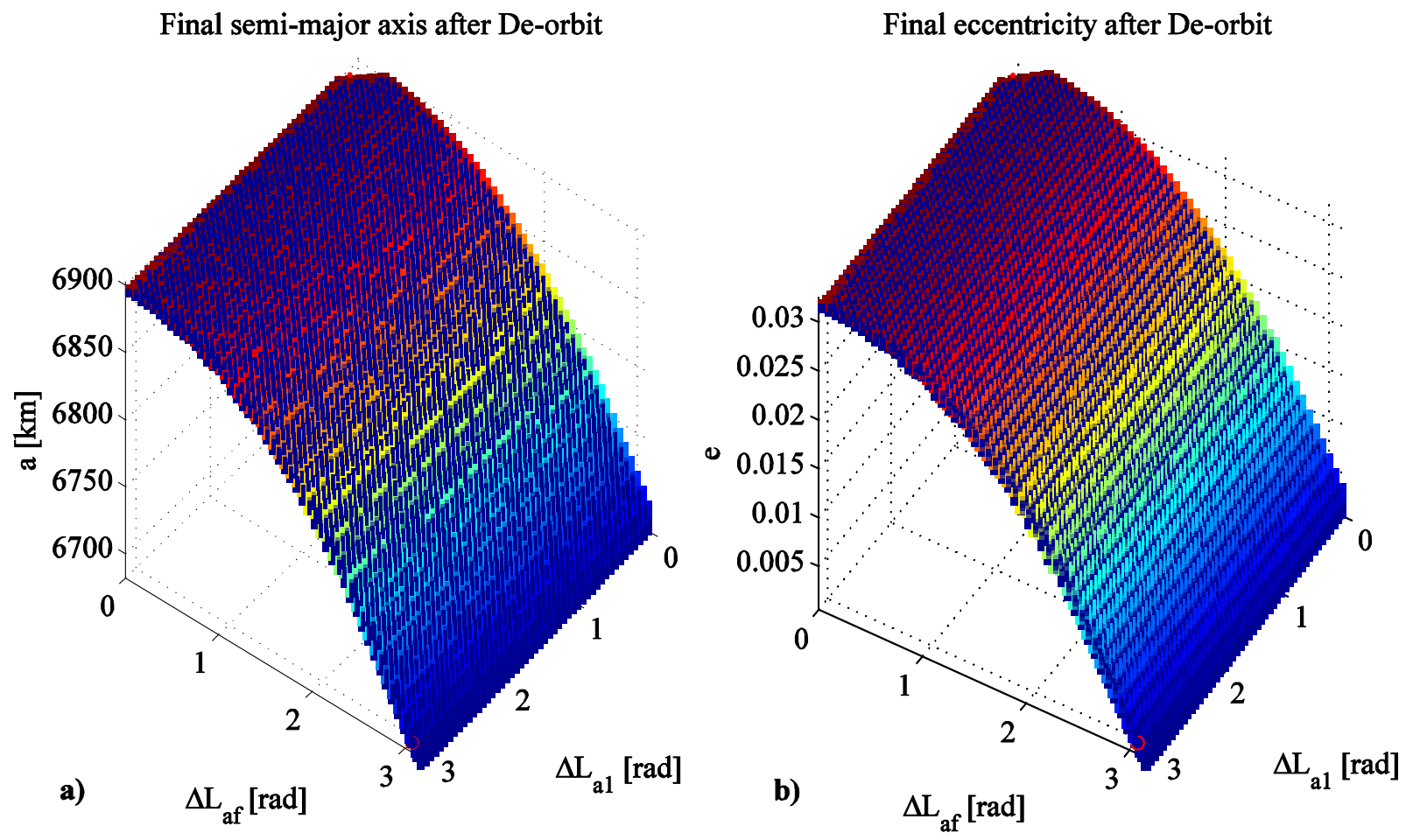

**Fig. 4 a)** final semi-major axis and **b)** eccentricity after de-orbit with respect to $\Delta L_{a1}$ and $\Delta L_{af}$ for $m_{IBS0} = 300kg$, $a_{debr0} = 7128km$ and $m_{debr} = 120kg$

Now it is desirable that the surrogate model returns the $\Delta V$ cost as a function of $m_{IBS0}$, $m_{debr}$, $a_{debr0}$ and *ToF*. From the available data relating the $\Delta V$ and *ToF* to the decision variables $\Delta L_{a1}$ and $\Delta L_{af}$ one can derive the functional relationship between $\Delta V$ and *ToF*. Given a triplet $m_{IBS0}$, $m_{debr}$, $a_{debr0}$, each *ToF* value defines a level curve on the $\Delta L_{a1}$ and $\Delta L_{af}$ plane (see Fig. 3a), which can be mapped into a set of $\Delta V$ values (see Fig.3b). Within this set, one can take the element with minimum $\Delta V$. Thus, for each time of flight, between a minimum and a maximum, one can derive the corresponding minimum $\Delta V$ cost. A similar procedure is followed to find the functional relationship between the final semi-major axis and the *ToF*. Note that there is no need to do the same for the eccentricity given the fact that the final perigee radius is fixed at $\overline{r}_p$ and therefore the final $e$ can be computed from the final $a$. In this way one can build the two surrogate models:

$$\begin{aligned} \Delta V &= f_{\Delta V,interp}\left(ToF, m_{IBS0}, m_{debr}, a_{debr0}\right) \\ a_f &= f_{a_f,interp}\left(ToF, m_{IBS0}, m_{debr}, a_{debr0}\right) \end{aligned} \tag{17}$$

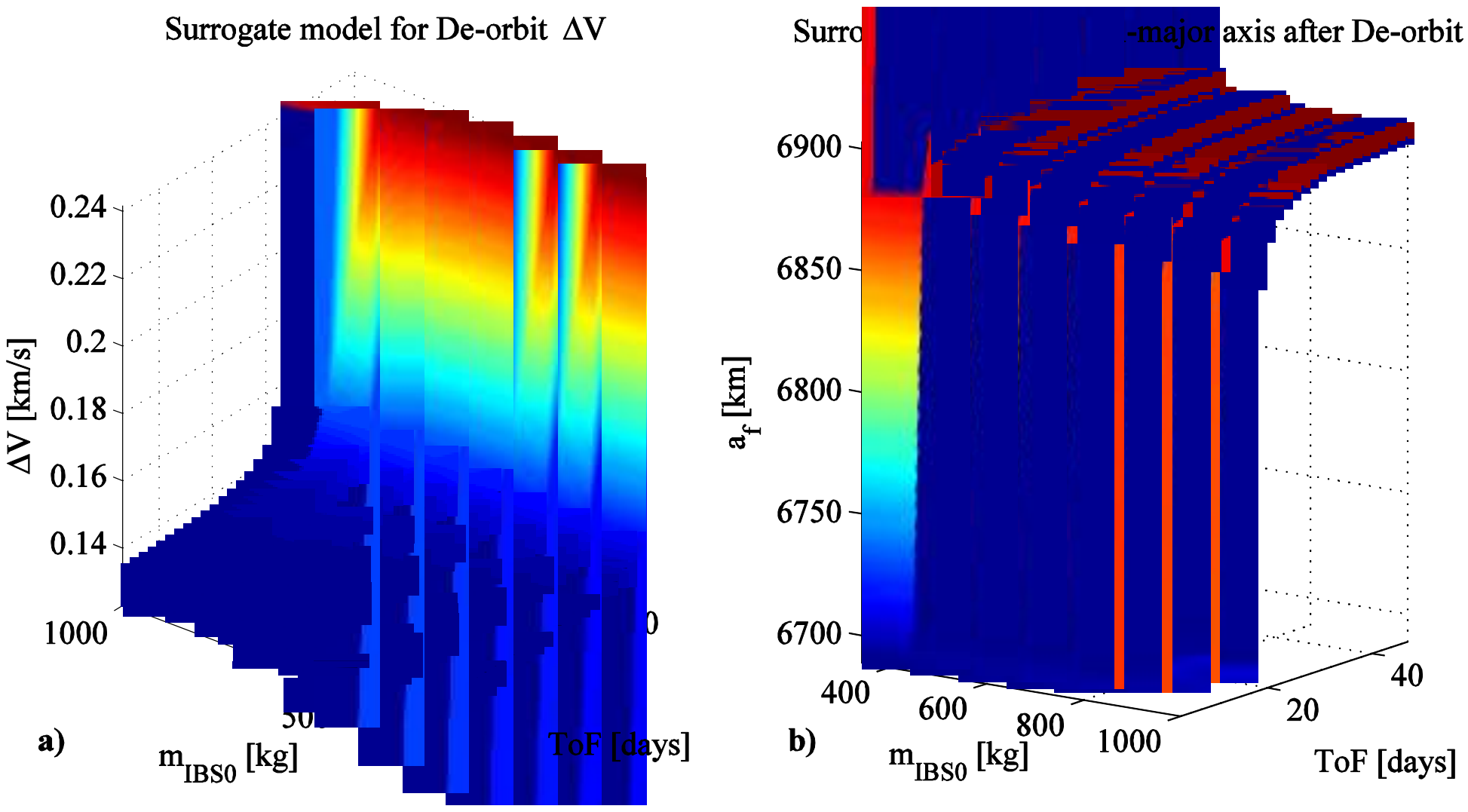

**Fig. 5** 3D Plot of surrogate models for $a_{debr0} = 7128km$ and $m_{debr} = 120kg$ : **a)** $\Delta V$; **b)** $a_f$

Fig. 5a and Fig. 5b show examples of tri-dimensional plots ($m_{IBS0}$-$ToF$-$\Delta V$ and $m_{IBS0}$-$ToF$-$a_f$ respectively) created by evaluating the surrogated models keeping $a_{debr0}$ and $m_{debr}$ fixed. In Fig. 5a one can see that there is a large *plateau* region corresponding to large time of flights and a smaller region close to the minimum *ToF* where the de-orbit cost increases very steeply and the final semi-major axis in Fig. 5b similarly decreases. The complete procedure for the creation of the interpolated de-orbit cost models requires few minutes of CPU time and once completed allows for a very fast estimation of the de-orbit cost. The surrogated models will be extremely useful in the multi-objective optimisation of debris removal sequences as it will be shown in the following sections.

### 3.2. Orbit Transfer Model

According to the scenario presented in Section 3, after having left the debris on a re-entry orbit, the IBS will have to transfer to the orbit of the next debris and rendezvous with it. The design of such a transfer arc would normally require the solution of a fixed-time Two Point Boundary Value Problem (TPBVP) which would be computationally very expensive given the high number of control parameters and constraints involved. A second simplified model was then created to quickly estimate the cost of a low-thrust multi-revolution orbit transfer with boundary constraints. The approach and assumptions presented in this section are similar to those already introduced for the de-orbit model.

First, given the limited acceleration provided by low thrust propulsion systems, one should consider that the orbit transfer will require a high number of multiple revolutions around the Earth, typical in the range of hundreds to few thousands. In this sense, it is possible to argue that achieving the proper phasing to transfer from the initial to final orbit would not be a major issue. Even a small variation of $\omega$ and $\theta$ per revolution would be sufficient to attain the required orientation to rendezvous with the piece of debris. Moreover, it is important to bear in mind that, in order to de-orbit the previous debris in the sequence, the IBS spacecraft, started from a circular orbit which was subsequently modified into an elliptical one with perigee $\overline{r}_p$. Thus it would be also possible to conveniently adjust the start point of the de-orbit procedure from the circular orbit in order to obtain the proper phasing once this is completed. For all these reasons, it is assumed that in this particular case, the phasing problem will have a negligible effect on the $\Delta V$ and time required to rendezvous with the next piece of debris in the sequence. Therefore, in the following it is assumed that it is not necessary to match the arrival $\omega$ and $\theta$ computed with the simplified model with those of the target object. Matching the target inclination $i$ and RAAN $\Omega$, instead, cannot be ignored without introducing a considerable error in the $\Delta V$ cost. In order to match the inclination and RAAN difference, one need to take into account only the geometric angle between planes of the initial and final orbits, which is given by:

$$\Delta i = \arccos\left(-\cos i_0 \cos\left(\pi - i_f\right) + \sin i_0 \sin\left(\pi - i_f\right)\cos\left(\Omega_f - \Omega_0\right)\right) \tag{18}$$

Thus in order to account for $\Delta i$, the inclination of the initial orbit is fictitiously set to zero, while the final one is set at $\Delta i$. The matching of the RAAN is assured by performing the circularisation properly. The assumption is that the deorbiting of one piece of debris starts at a true anomaly such that the resulting elliptical orbit has the line of apses perpendicular with the line of the nodes of the following piece of debris. Since the orbits of the debris are assumed to be circular, it is always possible to start the deorbiting at the right true anomaly with minimum delay. This hypothesis will be discussed in more detail with some numerical examples in Section 4.

With these assumptions, the main issue in designing the multi-revolution transfer will be that of achieving the required change in the apogee and perigee radiuses in order to match those of the final orbit, and to achieve the required rotation of the orbit plane.

The control of eccentricity and semi-major axis, required to match the target perigee and apogee altitudes, can be obtained by inserting two thrust arcs per revolution, one around the apogee and one around the perigee. This methodology is analogous to what was done in the previous section for the perigee lowering. In the same way, the radial component of the thrust acceleration is set to zero. The transverse component this time can have either positive or negative sign ( $\alpha = \pm\pi/2$ ) depending whether the perigee (or apogee) needs to be raised or lowered.

Since a plane change is required, the out-of-plane component of the thrust acceleration is non-zero. Thanks to this the control parameters can be reduced to the semi-amplitude of the apogee and perigee thrusting arcs, $\Delta L_a$ and $\Delta L_p$, the sign of the $\theta$ component of the thrust acceleration (i.e. the sign of $\alpha_a, \alpha_p = \pm\pi/2$ ) and the out-of-plane component in the same arcs, $\beta_a$ and $\beta_p$. Define $\Delta L_{thrust}$ as half the total thrusting arc length and $r_t$ as the ratio of $\Delta L_{thrust}$ which is devoted to apogee thrusting. In order to have a parameterisation which accounts also for the sign of $\alpha_a$ and $\alpha_p$ the following one is proposed:

$$
\begin{aligned}
\alpha_a &= \begin{cases} \pi/2 & \Delta L_t \geq 0 \\ -\pi/2 & \Delta L_t < 0 \end{cases} \\
\alpha_p &= \begin{cases} \alpha_a & 0 \leq r_t \leq 1 \\ -\alpha_a & 1 < r_t \leq 2 \end{cases} \\
\Delta L_a &= \begin{cases} r_t \left|\Delta L_{thrust}\right| & 0 \leq r_t \leq 1 \\ (2 - r_t)\left|\Delta L_{thrust}\right| & 1 < r_t \leq 2 \end{cases} \\
\Delta L_p &= \left|\Delta L_{thrust}\right| - \Delta L_a
\end{aligned}
\tag{19}
$$

with

$$
\begin{aligned}
\Delta L_{thrust} &\in [-\pi \quad \pi] \\
r_t &\in [0 \quad 2]
\end{aligned}
$$

To define the actual values of $\Delta L_{thrust}$ and $r_t$ in each revolution, an interpolating strategy from a set of nodal values similar to the one used for the de-orbit model is adopted. In this case, however, the interpolated values will not be computed with respect to the current revolution number but with respect to the time. Again the number of interpolating nodes can be chosen arbitrarily and is set to 2 in this case, $\Delta L_{t1}$, $\Delta L_{tf}$, $r_{t1}$, $r_{tf}$. For $\beta_a$ and $\beta_p$, it is chosen to have a constant value along the entire transfer. The thrusting pattern along each revolution is shown in Fig. 6.

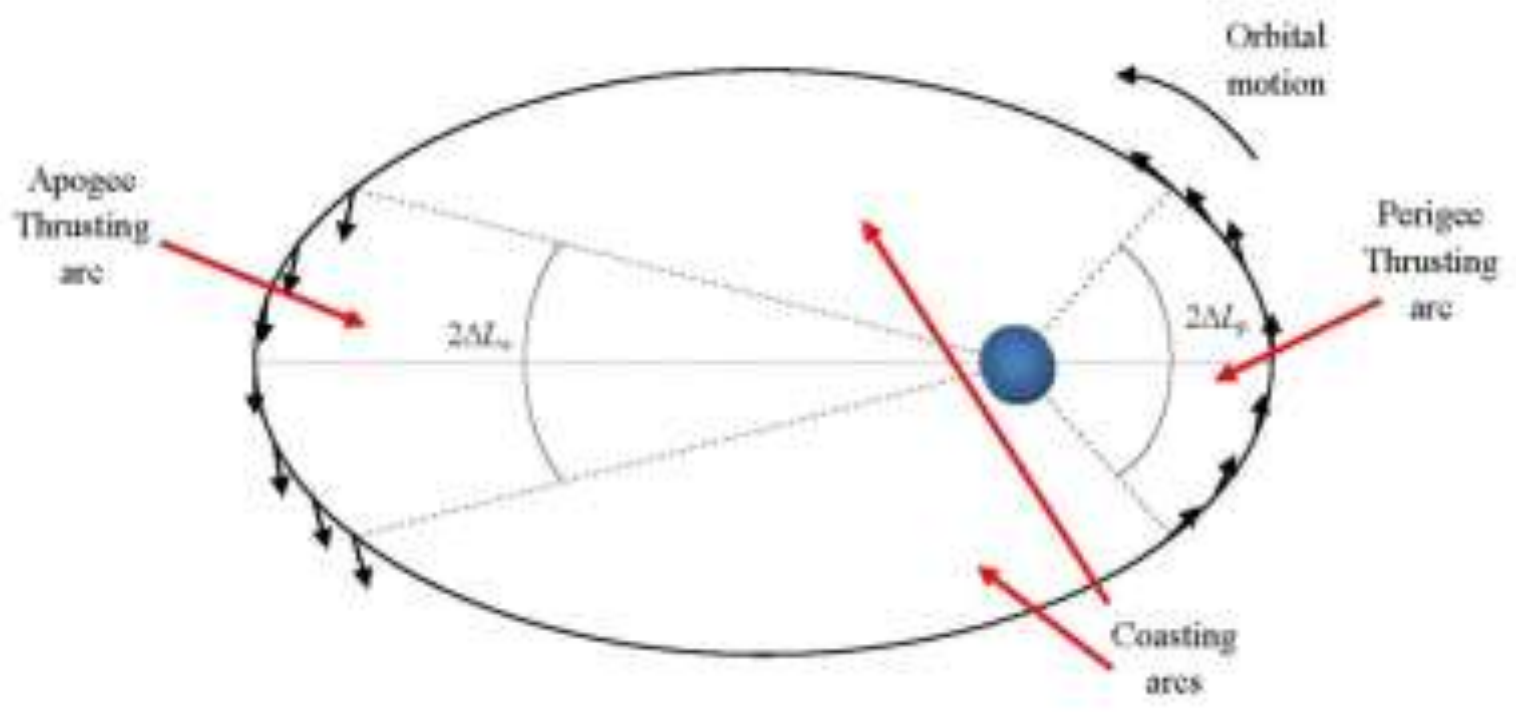

**Fig. 6** Thrusting arcs around apogee and perigee

Given a set of control parameters $\left[\Delta L_{t1} \quad \Delta L_{tf} \quad r_{t1} \quad r_{tf} \quad \beta_a \quad \beta_p\right]$ a multi-revolution transfer with specified duration $\overline{ToF}$, departing from an orbit defined by $\left[a_0 \quad e_0 \quad 0\right]^T$ and targeted to an orbit defined by $\left[a_f \quad e_f \quad \Delta i\right]^T$, is propagated according to the following procedure:

1. Compute the set of initial Equinoctial parameters $L_0$ and $\mathbf{E}_0 = \left[a_0 \quad P_{10} \quad P_{20} \quad Q_{10} \quad Q_{20}\right]^T$. $Q_{10}$ and $Q_{20}$ will be zero since the initial inclination is arbitrarily set to zero.
2. Compute the set of target Equinoctial parameters $\bar{\mathbf{E}}_f = \left[\bar{a}_f \quad \bar{P}_{1f} \quad \bar{P}_{2f} \quad \bar{Q}_{1f} \quad \bar{Q}_{2f}\right]^T$. Note that $\bar{P}_{1f}$ and $\bar{P}_{2f}$ will be zero since in this case the target orbit is a circular one.
3. Initialise the total $\Delta V$ and Time of flight to zero:

$$\Delta V = 0$$
$$ToF = 0$$

4. Set $\mathbf{E}_p^- = \mathbf{E}_0$ and $L_{coast,a} = L_0$.
5. Initialise the mass of the IBS spacecraft:

$$m_{IBS} = m_{IBS0}$$

6. While $ToF < \overline{ToF}$:
   a. Compute the interpolated values for $\Delta L_t$ and $r_t$. Hence calculate $\alpha_a$, $\alpha_p$, $\Delta L_a$ and $\Delta L_p$ from Eq. (19).
   b. Compute:

$$\begin{aligned} L_a^- &= L_a - \Delta L_a \qquad L_a^+ = L_a + \Delta L_a \\ L_p^- &= L_p - \Delta L_p \qquad L_p^+ = L_p + \Delta L_p \end{aligned} \tag{20}$$

   c. Compute the current acceleration acting on the spacecraft:

$$\varepsilon_{IBS} = \frac{F_{tot}}{m_{IBS}} \tag{21}$$

   d. Compute the time of flight $t_{coast,p}$ spent coasting before perigee from $L_{coast,p}$ to $L_p^-$.
   e. Compute the Equinoctial parameters after the thrusting perigee arc $\mathbf{E}_p^+$ with an expression analogous to Eq. (12).
   f. Compute the thrusting time at perigee $t_{thrust,p}$ from Eq. (13). If $\left(\overline{ToF} - ToF\right) < t_{thrust,p}$ proceed to step g. Otherwise, break the iterative sequence and go to step 7.
   g. Update $\Delta V$ and $ToF$:

$$\Delta V = \Delta V + \varepsilon_{IBS} t_{thrust,p} \tag{22}$$

$$ToF = ToF + t_{coast,p} + t_{thrust,p} \tag{23}$$

h. Update the IBS spacecraft mass:

$$m_{IBS} = m_{IBS} \exp\left( -\frac{\varepsilon_{IBS} t_{thrust,p}}{I_{sp} g_0} \right) \tag{24}$$

i. Set $\mathbf{E}_a^{-} = \mathbf{E}_p^{+}$ and $L_{coast,a} = L_p^{+}$.

j. Compute the current acceleration on the spacecraft:

$$\varepsilon_{IBS} = \frac{F_{tot}}{m_{IBS}} \tag{25}$$

k. Compute the time of flight $t_{coast,a}$ spent coasting before apogee from $L_{coast,a}$ to $L_a^{-}$.

l. Compute the Equinoctial parameters after the thrusting apogee arc $\mathbf{E}_a^{+}$ as in Eq. (12).

m. Compute the thrusting time at apogee $t_{thrust,a}$ from Eq. (13). If $\left(\overline{ToF} - ToF\right) < t_{thrust,a}$ proceed to step n. Otherwise, break the iterative sequence and go to step 7.

n. Update $\Delta V$ and $ToF$:

$$\Delta V = \Delta V + \varepsilon_{IBS} t_{thrust,a} \tag{26}$$

$$ToF = ToF + t_{coast,a} + t_{thrust,a} \tag{27}$$

o. Update the IBS spacecraft mass:

$$m_{IBS} = m_{IBS} \exp\left( -\frac{\varepsilon_{IBS} t_{thrust,a}}{I_{sp} g_0} \right) \tag{28}$$

p. Set $\mathbf{E}_p^{-} = \mathbf{E}_a^{+}$ and $L_{coast,p} = L_a^{+}$.

7. Back-track the point at which $ToF = \overline{ToF}$ and compute the corresponding equinoctial parameters $\mathbf{E}_f = \left[ a_f \quad P_{1f} \quad P_{2f} \quad Q_{1f} \quad Q_{2f} \right]^T$ and update $\Delta V$ accordingly.
8. Compute the mismatch between the actual final conditions and the target orbit:

$$C_{eq} = \begin{bmatrix} a_f - \bar{a}_f \\ e_f - \bar{e}_f = \sqrt{P_{1f}^2 + P_{2f}^2} - \sqrt{\bar{P}_{1f}^2 + \bar{P}_{2f}^2} \\ i_f - \bar{i}_f = 2\left( \arctan\sqrt{Q_{1f}^2 + Q_{2f}^2} - \arctan\sqrt{\bar{Q}_{1f}^2 + \bar{Q}_{2f}^2} \right) \end{bmatrix} \tag{29}$$

Summarizing, the TPBVP has been reduced to an optimisation problem in the form:

$$\begin{gathered} \min_{\mathbf{x}} \Delta V \\ s.t.\, C_{eq} = 0 \\ with\ \mathbf{x} = \left[ \Delta L_{t1}, \Delta L_{tf}, r_{t1}, r_{tf}, \beta_a, \beta_p \right] \end{gathered} \tag{30}$$

This problem can be solved with a gradient-based optimisation algorithm like MatLab®'s *fmincon*. Note that, the time of flight $\overline{ToF}$ is specified a priori and therefore it might occur that this duration is too short as to obtain the change in the orbital parameters specified by the boundary constraints. In this case, the problem is infeasible and the optimisation is terminated after a maximum 50 if the constraints are not satisfied.

In the following, an example of transfer from an elliptical orbit with 300 km perigee altitude and eccentricity 0.031 (corresponding to the final orbit of a de-orbiting strategy) to a circular orbit of 1100 km altitude (corresponding to the orbit of the next debris in an hypothetical removal sequence). Parameters of the two orbits are reported in Table 2. Note that the total plane rotation Δ*i* in this case is 10 degrees. The specified time of flight is 70 days.

**Table 2** Parameters of departure and arrival orbits

| | *a* [km] | *e* | *i* [deg] |
|---|---|---|---|
| Departure | 6892.24 | 0.031 | 0 |
| Arrival | 7478.16 | 0 | 10 |

First it is considered the case of a coplanar transfer, i.e. Δ*i*=0 will be computed. The optimisation problem was solved with *fmincon* in 6 iterations and less than 10 seconds, returning a minimum Δ*V* cost of 0.301 km/s, with 1001 revolutions. Fig. 7a-c report respectively the variation of semi-major axis, eccentricity, apogee and perigee radii. One can see that *a* is monotonically increasing while *e* on the other hand is monotonically decreasing to zero. In order to reach the desired circular orbit, the perigee had to be raised by almost 700 km while the apogee had to be raised by some 400 km. This higher effort needed to raise the perigee explains the larger amplitude of apogee thrusting arcs $\Delta L_a$ compared to perigee ones $\Delta L_p$ (as shown in Fig. 8a). The azimuth thrust angles $\alpha_p$, $\alpha_a$ (see Fig. 8b) are both positive since both the perigee and apogee are raised. $\beta_p$ and $\beta_a$ are obviously zero because the transfer is coplanar and thus Δ*i* is constantly null.

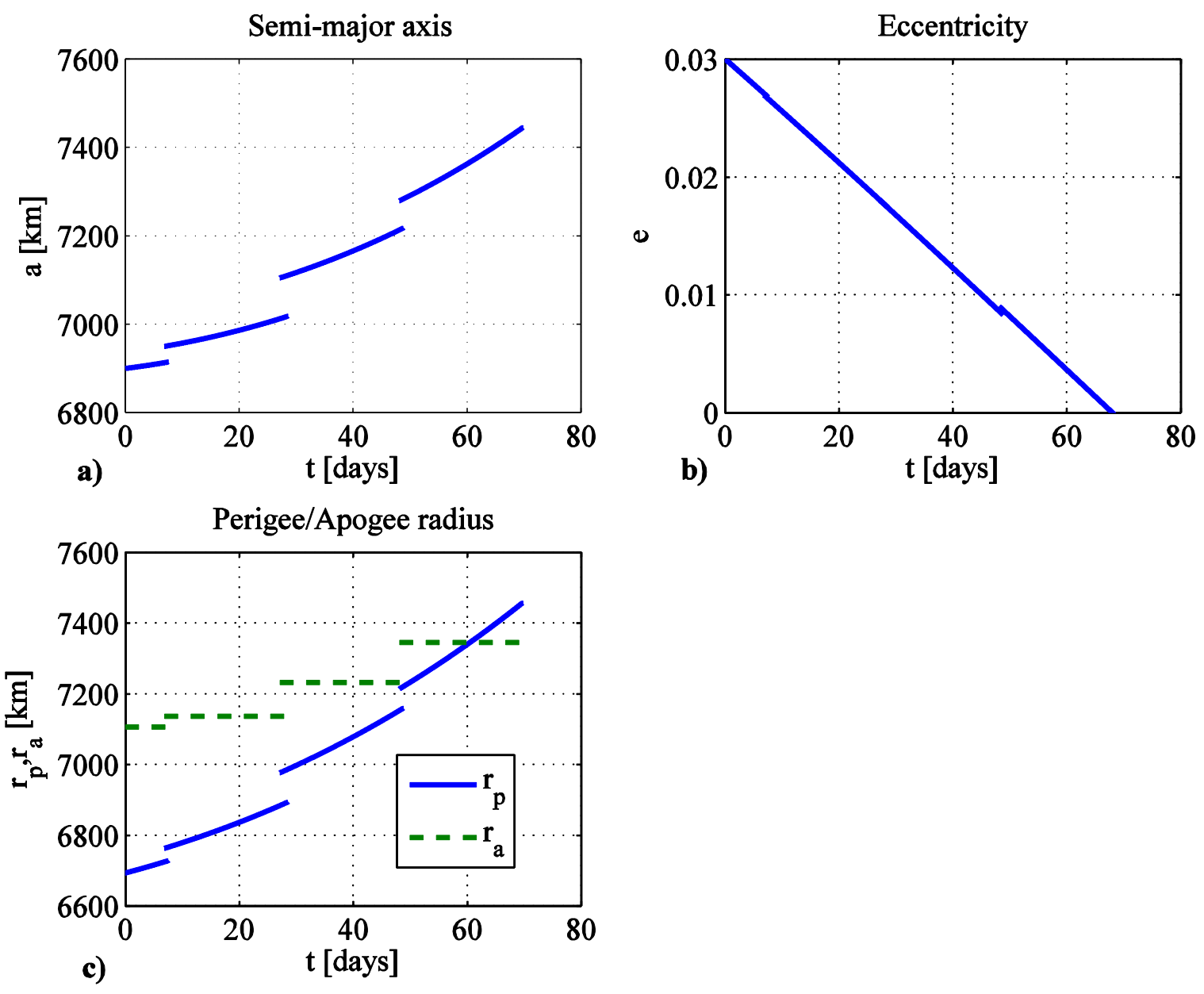

**Fig. 7** **a)** variation of semi-major axis, **b)** eccentricity, **c)** perigee and apogee radiuses for multi-revolution orbital transfer (coplanar case)

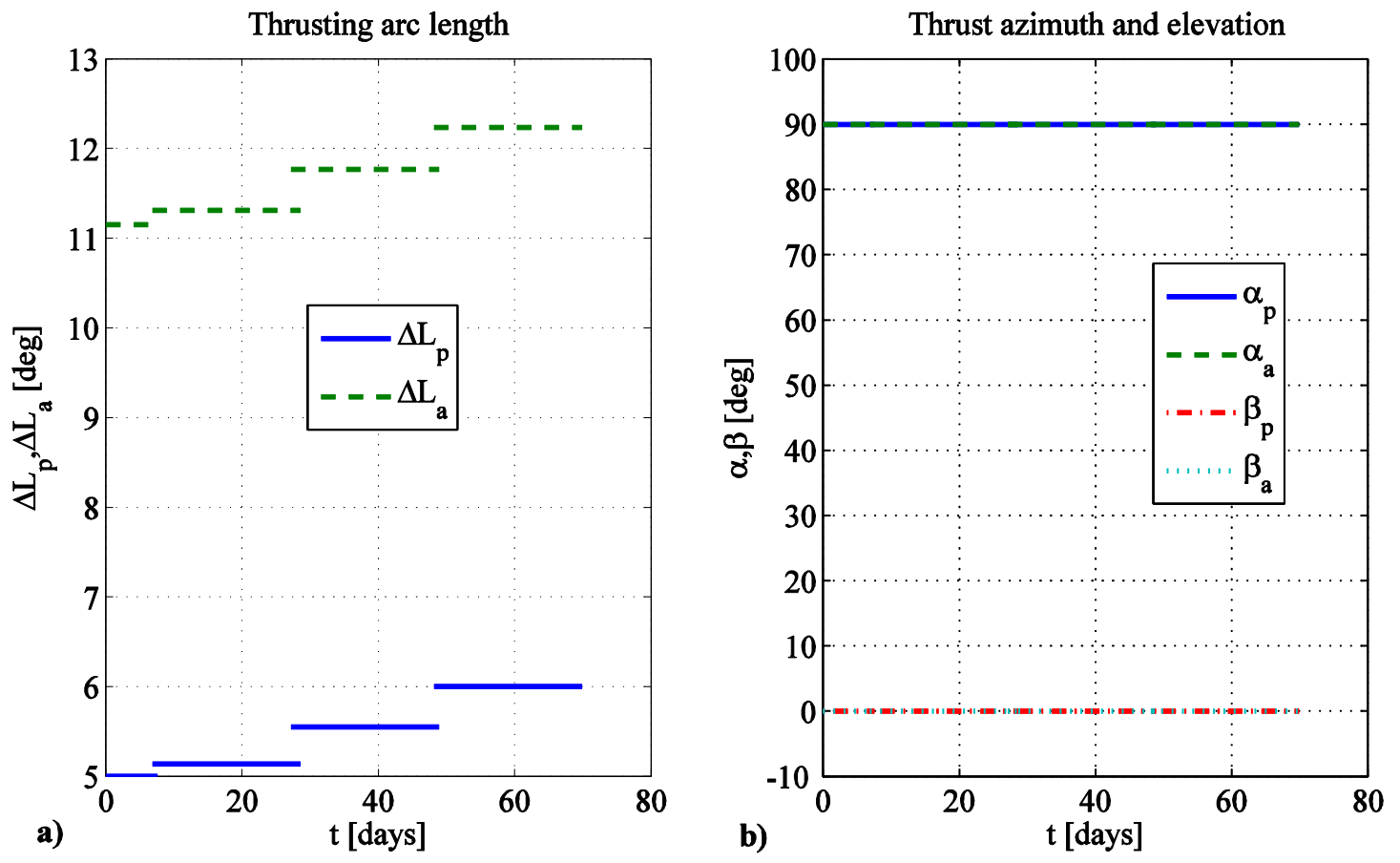

**Fig. 8** Control parameters for multi-revolution orbital transfer (coplanar case): **a)** thrust arc length; **b)** azimuth and elevation

The same problem, but this time with the 10° plane change specified in Table 2 returns a $\Delta V$ of 1.480 km/s with 1004 revolutions. The high cost of out plane manoeuvres is well exemplified by the fact that the $\Delta V$ required is more than four times larger than a coplanar transfer. As can be seen in Fig. 9a, Fig. 9b, Fig. 9d, semi-major axis, eccentricity, apogee and perigee radii show a similar behaviour to the coplanar case while this time also the inclination (as in Fig. 9c) increases monotonically to 10 degrees. By analysing the control parameters in Fig. 10a one can see that this time the amplitude of the perigee arcs in general larger than the apogee ones, even if, like in the coplanar case, the increase in perigee is much larger than that of the apogee. This  fact is explained by the fact that the out-of-plane component at perigee $\beta_p$ is close to 90° (see Fig. 10b), meaning that the thrusting action at perigee is mostly devoted to the plane change. In contrast, $\beta_a$ is smaller in magnitude, around -70° (the opposite sign is due to the fact that it is advantageous to invert the out-of-plane component twice per revolution), therefore with a higher in plane component devoted to perigee raising.

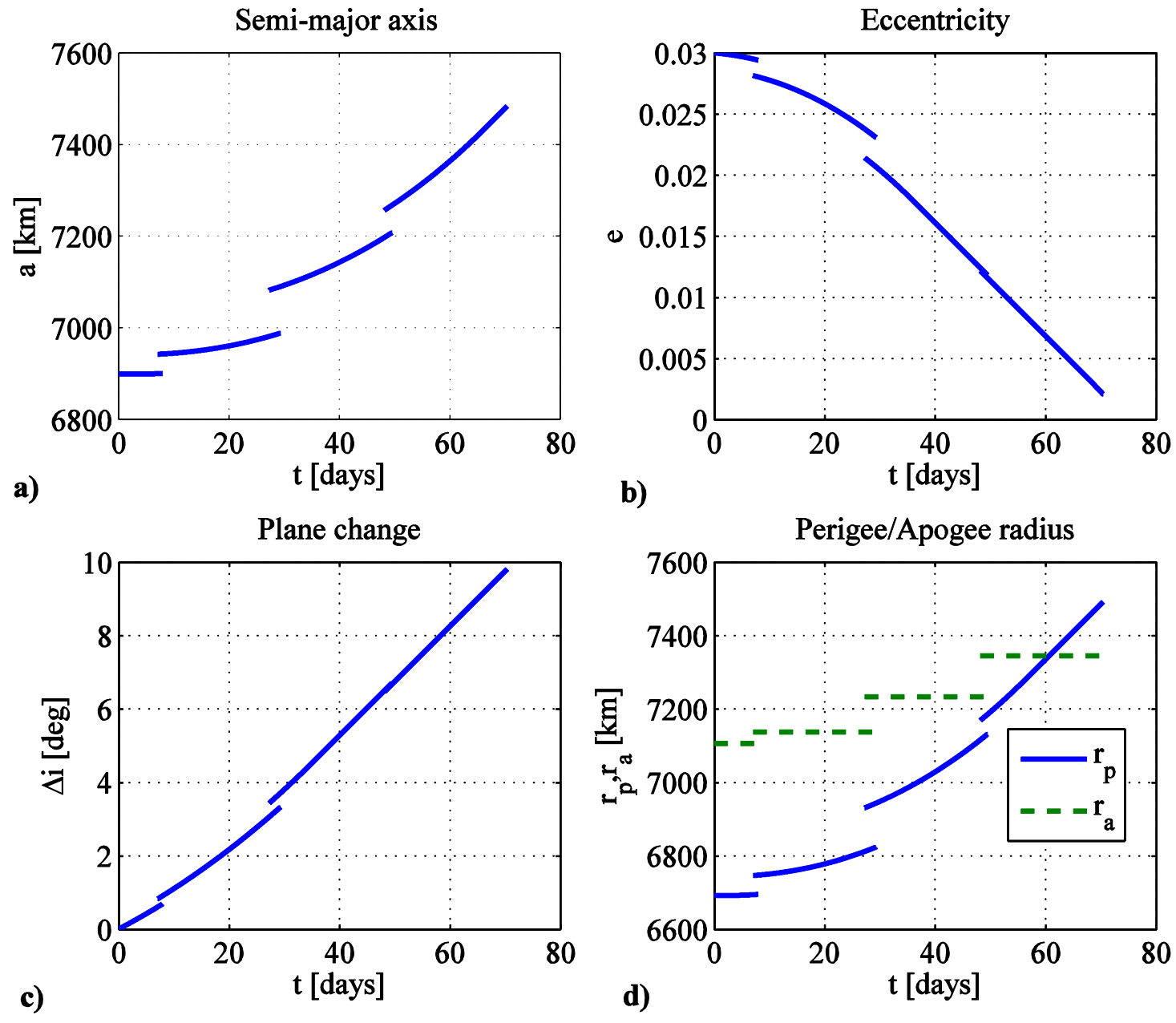

**Fig. 9 a)** variation of semi-major axis, **b)** eccentricity, **c)** plane change, **d)** perigee and apogee radiuses for multi-revolution orbital transfer (10° plane change)

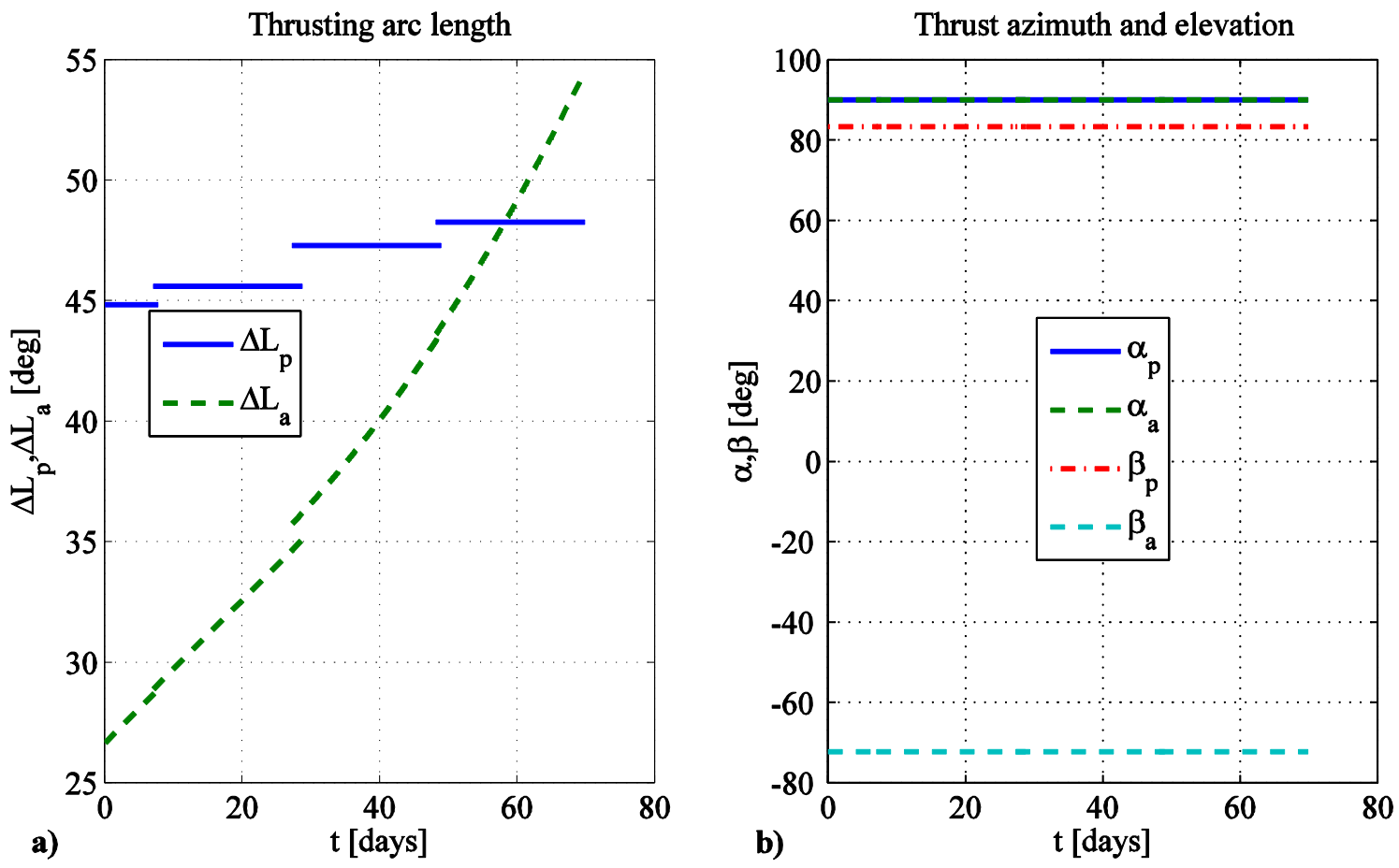

**Fig. 10** Control parameters for multi-revolution orbital transfer (10º plane change): **a)** thrust arc length; **b)** azimuth and elevation

## 4. Multi-Objective optimisation

The aim is now that of optimising the timing and sequence of a removal mission by means of a single IBS spacecraft. It is assumed that the spacecraft departs from a LEO with a 250 km semi-major axis altitude and coplanar with respect to the first piece of debris in the sequence. The five target objects have the orbital parameters and mass reported in Table 3. The mass and orbital parameters have been chosen arbitrarily while adhering to the observations in [9] and [11] that the most dangerous debris are located in LEO and generally weigh a few hundred kilos. Different values for $i$ and $\Omega$ are also taken in order to consider the fact that the pieces of debris, in principle, will be orbiting on different planes. Note that $T_{DO,min}$ has been computed with the procedure detailed in Section 3.1 and therefore depends on the characteristics of the IBS spacecraft. Moreover, it is also important to remark that these are only *best case* figures values which were computed with a *minimum* hypothetical wet mass of 350 kg (much lower than the actual launch mass of 1000 kg). The surrogate models in Eqs. (17) can in general consider wet masses between 350 kg and 1000 kg, as shown for example in Fig. 5a-b.

**Table 3** Mass, initial orbit parameters and minimum de-orbit time of the debris

| Debris nr. | *mass* [kg] | *a* [km] | *e* | *i* [deg] | *Ω* [deg] | $T_{DO,min}$ [days] |
|---|---|---|---|---|---|---|
| 1 | 500 | 6828.16 | 0 | 1 | 65 | 2.67 |
| 2 | 120 | 7128.16 | 0 | 2 | 150 | 3.36 |
| 3 | 300 | 6978.16 | 0 | -2 | 200 | 3.68 |
| 4 | 400 | 7478.16 | 0 | -1 | 90 | 11.12 |
| 5 | 800 | 7178.16 | 0 | 0 | 45 | 12.25 |

Table 4 reports the relative inclination change between the orbit planes of the 5 different objects, as computed from Eq. (18).

**Table 4** Relative inclination change |Δ*i*| [deg] between orbit planes of the debris

| Debris nr. | 2 | 3 | 4 | 5 |
|---|---|---|---|---|
| 1 | 2.16 | 1.47 | 1.95 | 1 |
| 2 | - | 3.63 | 2.65 | 2 |
| 3 | - | - | 2.52 | 2 |
| 4 | - | - | - | 1 |

The de-orbit sequence is defined by the order according to which the five pieces of debris are removed, the time needed to rendezvous with $T_{RV}$ and the time to de-orbit $T_{DO}$ each of them. The order is defined by the integer vector:

$$\mathbf{ord} = \begin{bmatrix} i_1 & i_2 & i_3 & i_4 & i_5 \end{bmatrix} \tag{31}$$

which collects the indexes of the objects in the a single debris removal sequence. Since there are five objects, there are 120 possible de-orbit sequences. The other parameters are contained in the vector **x**:

$$\mathbf{x}=\begin{bmatrix} T_{RV,i_1} & T_{DO,i_1} & T_{RV,i_2} & T_{DO,i_2} & T_{RV,i_3} & T_{DO,i_3} & T_{RV,i_4} & T_{DO,i_4} & T_{RV,i_5} & T_{DO,i_5} \end{bmatrix} \tag{32}$$

The performance of each sequence is assessed according to its total $\Delta V_{Tot}$ cost and time of flight $ToF_{Tot}$. The latter is computed simply as:

$$ToF_{Tot}=\sum \mathbf{x} \tag{33}$$

The total $\Delta V$ cost is calculated sequentially by adding up the costs of each of the ten phases (rendezvous and de-orbit for each debris). In particular, the cost of the rendezvous $\Delta V_{RV}$ is computed by solving the optimisation problem (30) and the de-orbit cost $\Delta V_{DO}$ is calculated from the surrogated model in Eq. (17). The final conditions after de-orbit are also computed from Eq. (17) since they will be the departure conditions for the following rendezvous step. The propellant mass consumption is also taken into account and updated throughout the entire sequence computation. In order to have only a real valued optimisation problem, it is chosen here to treat each of the 120 sequences as a bi-objective optimisation problem with **ord** fixed and ten design variables defined in **x**. Therefore, optimisation problem becomes:

$$\min_{\mathbf{x}\in D}\left[ToF_{Tot}(\mathbf{x}) \quad \Delta V_{Tot}(\mathbf{x})\right] \tag{34}$$

The domain $D$ is defined by the upper and lower boundaries defined in Table 5. Note that the lower boundaries for de-orbit time are set according to the sequence and the minimum times reported in Table 3.

**Table 5** Optimisation boundaries

| Parameter | $T_{RV,i1}$ | $T_{DO,i1}$ | $T_{RV,i2}$ | $T_{DO,i2}$ | $T_{RV,i3}$ | $T_{DO,i3}$ | $T_{RV,i4}$ | $T_{DO,i4}$ | $T_{RV,i5}$ | $T_{DO,i5}$ |
|---|---|---|---|---|---|---|---|---|---|---|
| Lower Bound | 5 | $T_{DO,min,i1}$ | 5 | $T_{DO,min,i2}$ | 5 | $T_{DO,min,i3}$ | 5 | $T_{DO,min,i4}$ | 5 | $T_{DO,min,i5}$ |
| Upper Bound | 100 | 50 | 100 | 50 | 100 | 50 | 100 | 50 | 100 | 50 |

Each bi-objective optimisation problem is solved with MACS, a hybrid-memetic stochastic optimisation algorithm [21]. MACS was run for 40000 function evaluations with 30 agents. Each of the 120 optimisation instances required roughly 6 days of computational time to complete. The outputs are represented by the Pareto optimal solutions w.r.t. $\Delta V_{Tot}$ and $ToF_{Tot}$. Fig. 11 to Fig. 15 collect the Pareto fronts according to the number of the first object in the sequence, i.e. the first index in the vector **ord**, as introduced in Eq. (31). In each figure, each colour represents the Pareto front corresponding to one of the 24 debris removal sequences starting with the same object. For example, Fig. 11 includes the Pareto fronts of sequences 12345, 13245, 14235, 15234, 12435 etc. .

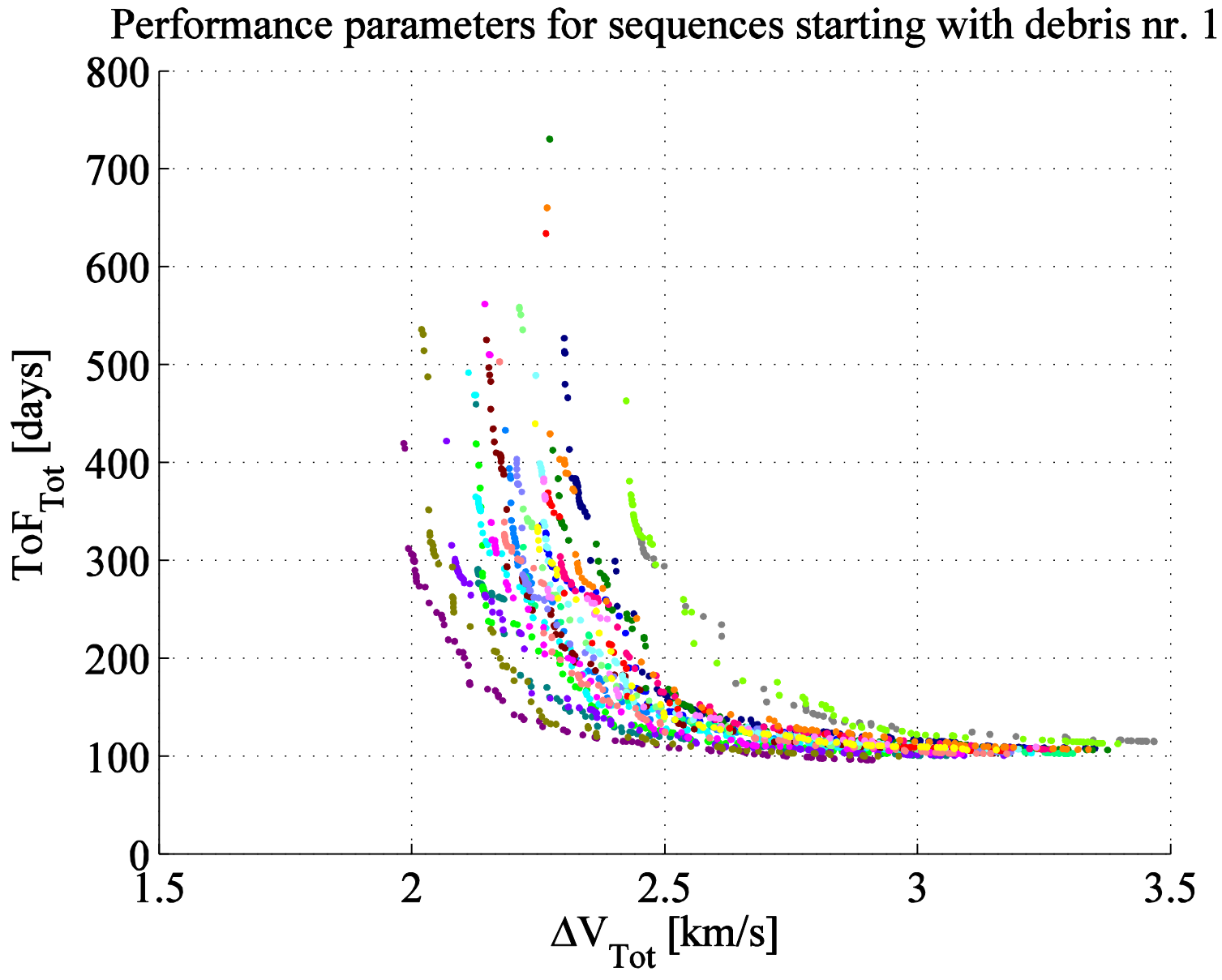

**Fig. 11** Pareto fronts for sequences starting with debris nr.1

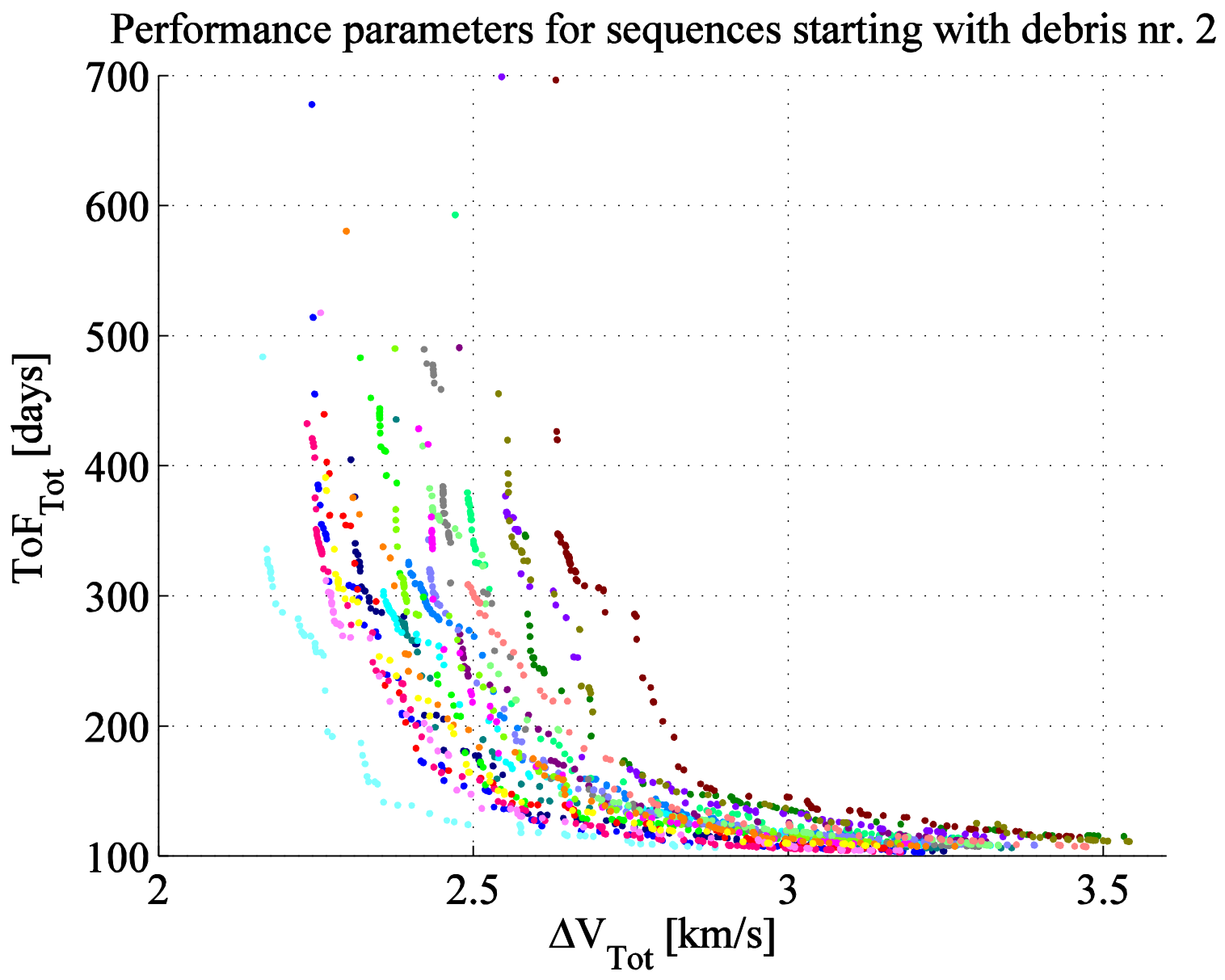

**Fig. 12** Pareto fronts for sequences starting with debris nr.2

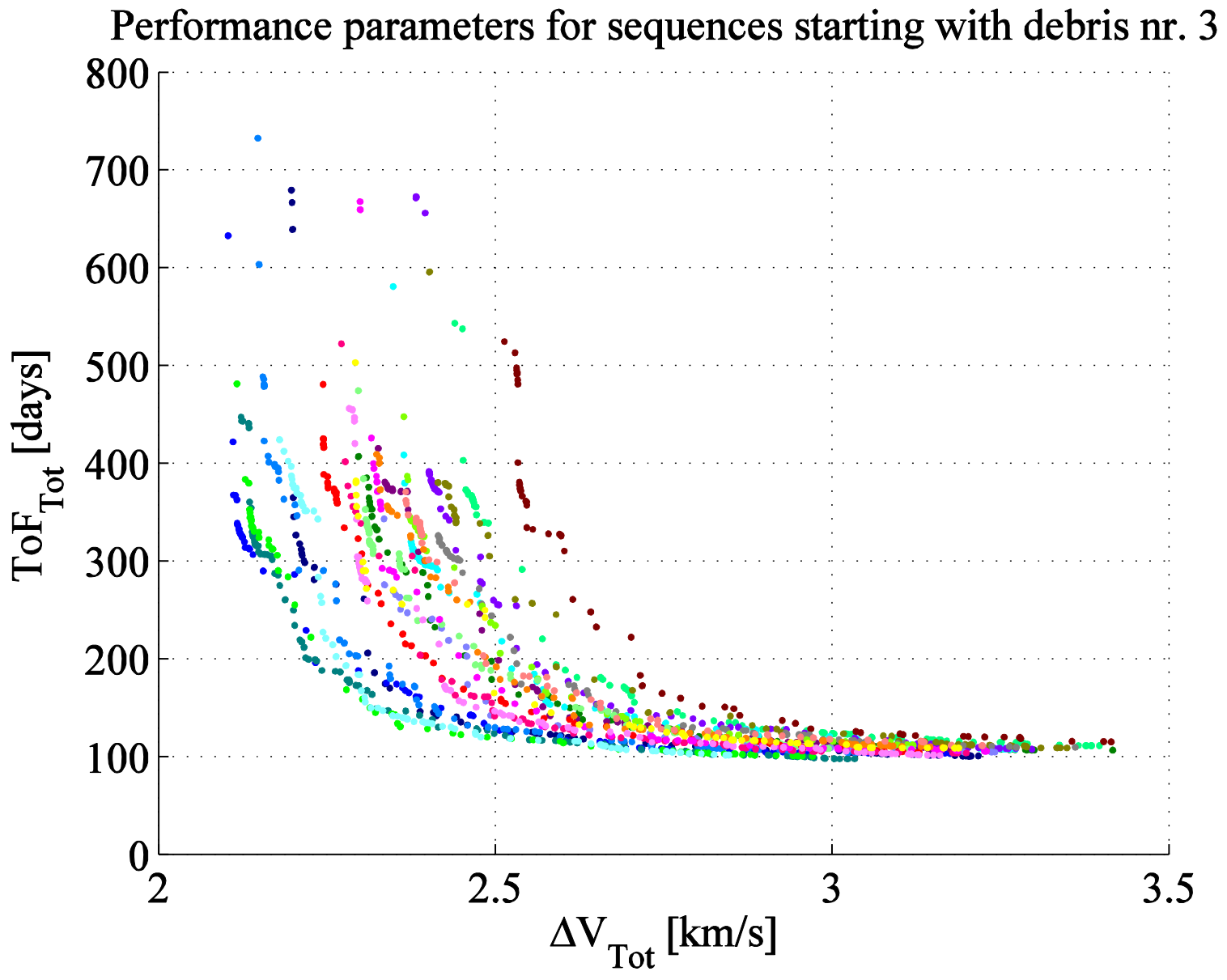

**Fig. 13** Pareto fronts for sequences starting with debris nr. 3

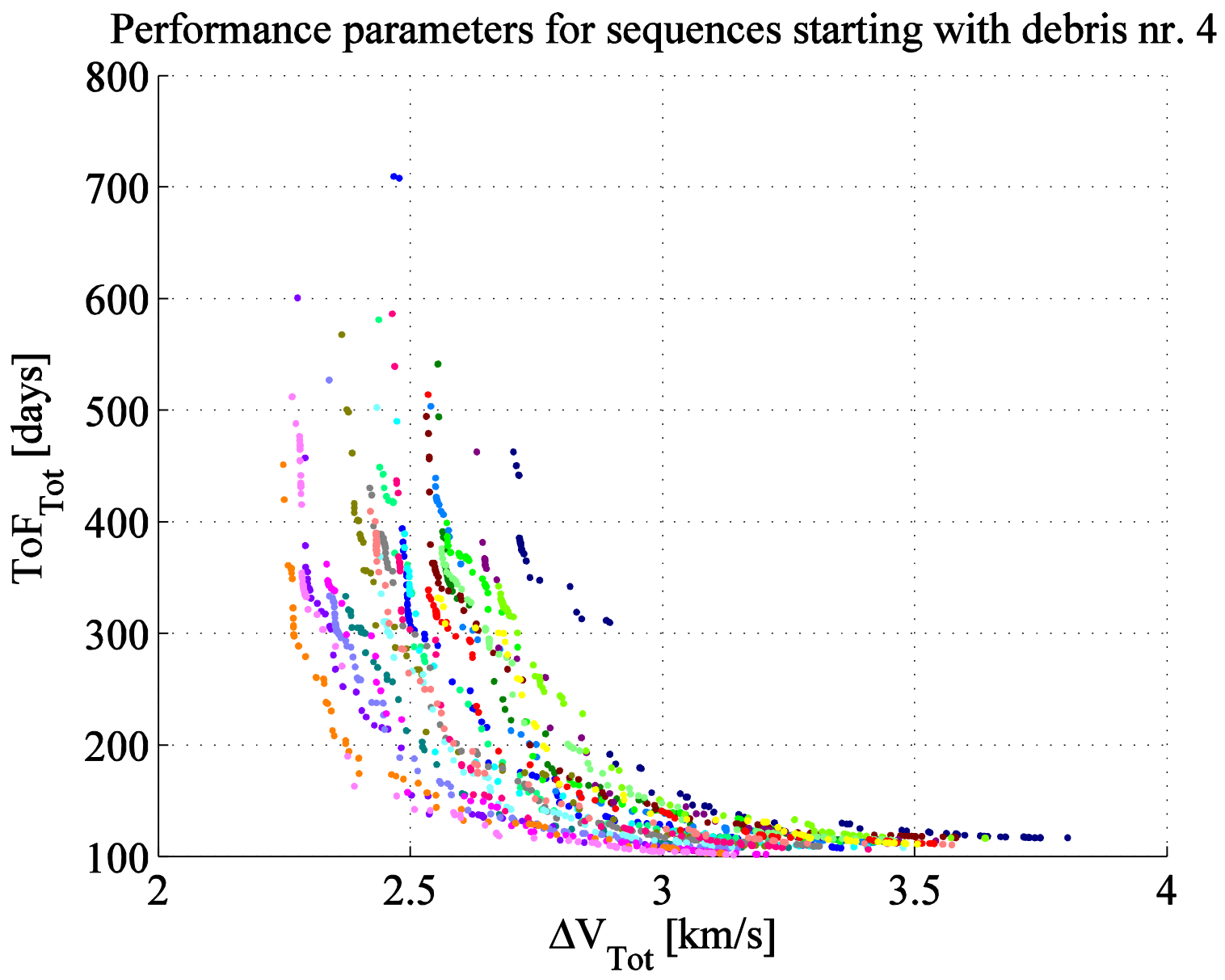

**Fig. 14** Pareto fronts for sequences starting with debris nr. 4

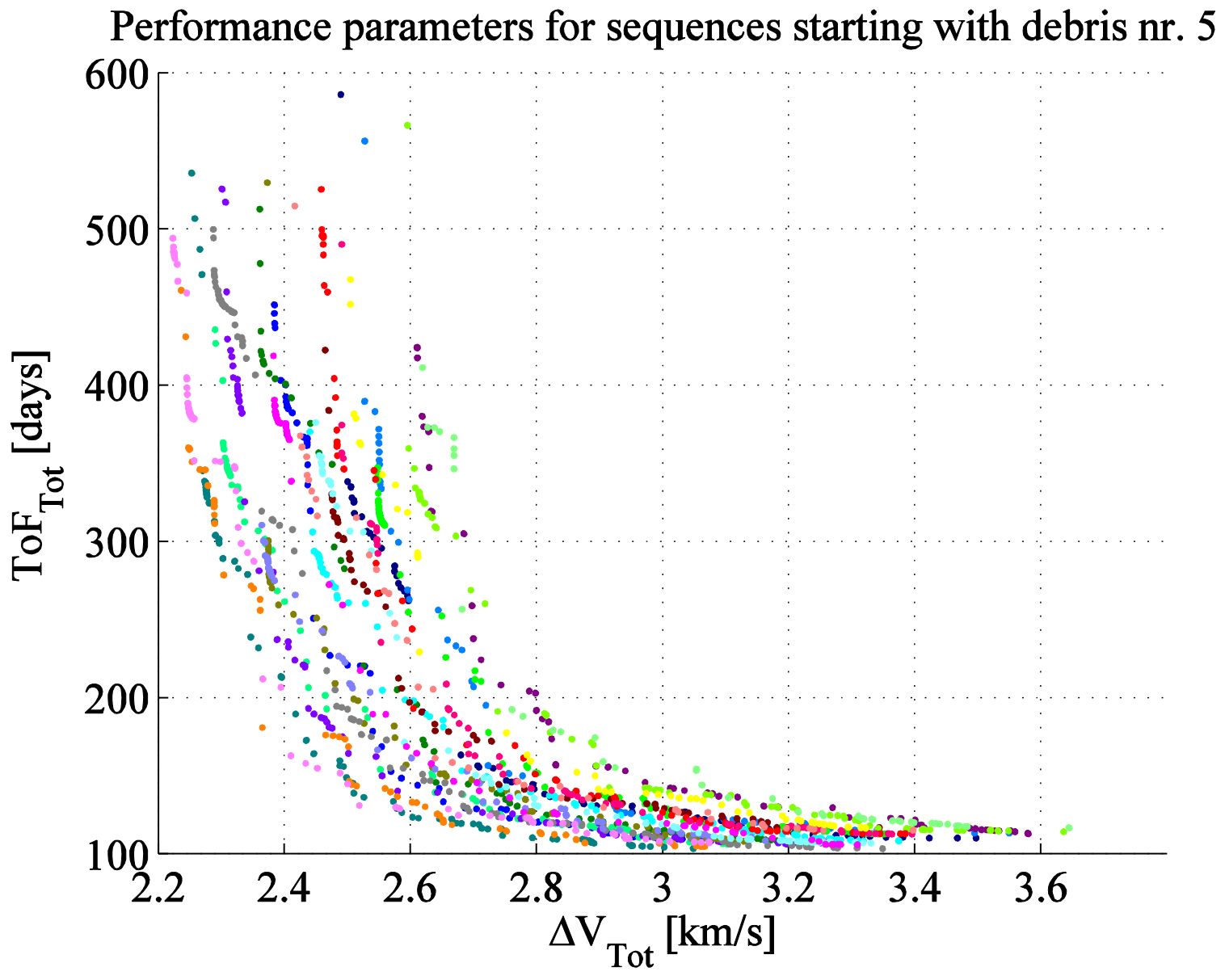

**Fig. 15** Pareto fronts for sequences starting with debris nr. 5

From a visual inspection of the fronts it is possible to see that sequences starting from debris nr. 1 seem to present the best $\Delta V_{Tot}$-$ToF_{Tot}$ combination, since for most of them the $\Delta V$ cost is comprised between 2 and 2.5 km/s. The corresponding times of flight are comprised roughly between 100 and 500 days. The sequences starting with debris nr. 3 and nr. 2 also have a good $\Delta V$ while those starting with nr. 4 and nr. 5 appear to be worst. By combining all the partial Pareto fronts one obtains the globally optimal solutions, as reported in Fig. 16.

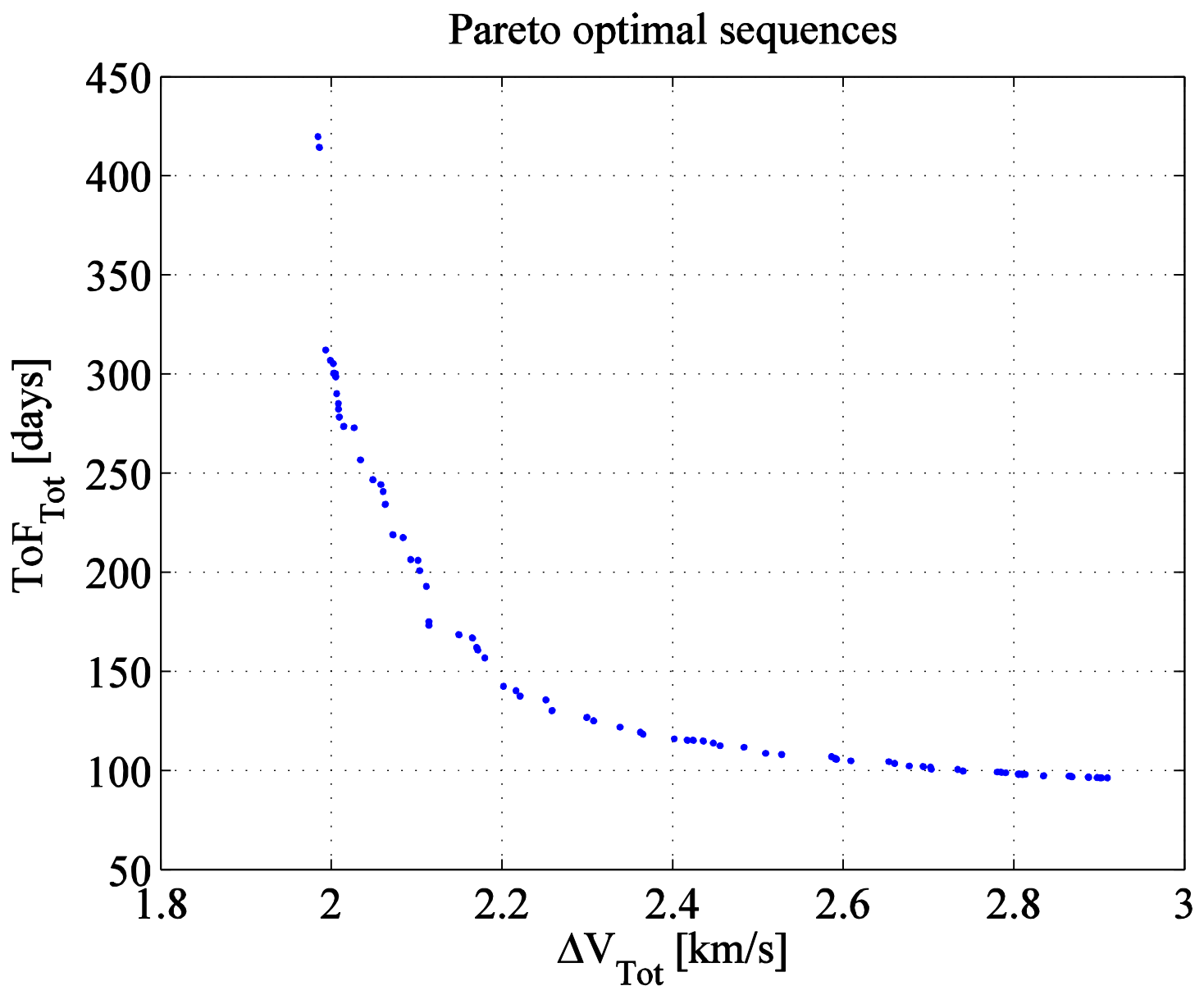

**Fig. 16** Global Pareto front

One can see that the global Pareto front is composed by individual solutions belonging exclusively from sequence 13452, which is therefore globally dominant. In order to rank the degree of optimality of each sequence, it is proposed to use an approach inspired by the performance metrics for optimisation algorithms proposed in [22]. Define $PF_g$ as the set of the points of the globally optimal Pareto front while $PF_{\mathbf{ord}}$ is the set of points belonging to the Pareto front corresponding to sequence **ord**. Define then the ranking parameter of sequence **ord** as:

$$Conv(\mathbf{ord}) = \frac{100\sum_{i=1}^{N_{\mathbf{ord}}} \min_{\mathbf{g}_j \in PF_g}\left(\frac{\mathbf{f}_i - \mathbf{g}_j}{\boldsymbol{\delta}}\right)}{N_{\mathbf{ord}}} \quad (35)$$
$$\mathbf{f}_i \in PF_{\mathbf{ord}}, \forall i = 1,\ldots,N_{\mathbf{ord}}$$

*Conv* is given by averaging the distance of each point of $PF_{\mathbf{ord}}$ from the closest point of $PF_g$. The closest $PF_{\mathbf{ord}}$ is to $PF_g$ and the lower *Conv* will be. Table 6 reports the ranking of the sequences according to *Conv*.

**Table 6** Ranking of the de-orbit sequences

| Rank | **ord** | *Conv*(**ord**) | Rank | **ord** | *Conv*(**ord**) | Rank | **ord** | *Conv*(**ord**) |
|---|---|---|---|---|---|---|---|---|
| 1 | 13452 | 0 | 41 | 42513 | 21.18 | 81 | 52143 | 31.43 |
| 2 | 13542 | 5.14 | 42 | 15234 | 21.26 | 82 | 32145 | 31.62 |
| 3 | 13524 | 6.61 | 43 | 32451 | 21.46 | 83 | 54123 | 31.72 |
| 4 | 12453 | 6.78 | 44 | 52134 | 21.46 | 84 | 54132 | 31.83 |
| 5 | 12543 | 7.25 | 45 | 34521 | 21.52 | 85 | 42135 | 32.43 |
| 6 | 31542 | 9.41 | 46 | 35142 | 21.79 | 86 | 52314 | 32.70 |
| 7 | 31452 | 9.85 | 47 | 35214 | 21.99 | 87 | 42531 | 33.05 |
| 8 | 34512 | 11.59 | 48 | 34251 | 22.02 | 88 | 21435 | 33.96 |
| 9 | 24513 | 12.15 | 49 | 52431 | 22.04 | 89 | 54231 | 34.05 |
| 10 | 15243 | 12.16 | 50 | 45132 | 22.13 | 90 | 23145 | 34.31 |
| 11 | 12534 | 12.33 | 51 | 54312 | 23.39 | 91 | 23514 | 34.56 |
| 12 | 31254 | 12.37 | 52 | 21543 | 23.60 | 92 | 53421 | 34.67 |
| 13 | 15432 | 13.24 | 53 | 24315 | 23.62 | 93 | 25341 | 34.71 |
| 14 | 35124 | 13.87 | 54 | 41352 | 23.81 | 94 | 14325 | 34.91 |
| 15 | 13254 | 14.22 | 55 | 43152 | 23.90 | 95 | 41253 | 35.20 |
| 16 | 31524 | 14.36 | 56 | 12435 | 24.40 | 96 | 32514 | 35.42 |
| 17 | 15342 | 14.48 | 57 | 34125 | 24.53 | 97 | 14235 | 35.65 |
| 18 | 13425 | 16.30 | 58 | 15324 | 24.89 | 98 | 32541 | 36.42 |
| 19 | 24531 | 16.53 | 59 | 53142 | 24.90 | 99 | 51234 | 36.81 |
| 20 | 14523 | 16.65 | 60 | 23154 | 25.61 | 100 | 42153 | 36.91 |
| 21 | 14352 | 16.69 | 61 | 53124 | 25.67 | 101 | 51423 | 38.10 |
| 22 | 34152 | 17.16 | 62 | 51243 | 25.80 | 102 | 54321 | 38.25 |
| 23 | 25134 | 17.17 | 63 | 43512 | 25.83 | 103 | 45231 | 40.16 |
| 24 | 12354 | 17.47 | 64 | 31425 | 25.95 | 104 | 51432 | 40.98 |
| 25 | 14253 | 17.63 | 65 | 12345 | 25.96 | 105 | 41523 | 41.91 |
| 26 | 31245 | 17.81 | 66 | 21453 | 26.01 | 106 | 45321 | 44.72 |
| 27 | 15423 | 17.85 | 67 | 52413 | 26.09 | 107 | 32415 | 45.05 |
| 28 | 51342 | 17.88 | 68 | 51324 | 26.56 | 108 | 42351 | 45.38 |
| 29 | 14532 | 18.05 | 69 | 35241 | 26.68 | 109 | 43521 | 45.43 |
| 30 | 25413 | 18.07 | 70 | 25143 | 26.77 | 110 | 53214 | 45.72 |
| 31 | 54213 | 18.17 | 71 | 24153 | 26.93 | 111 | 43251 | 45.87 |
| 32 | 35412 | 18.48 | 72 | 34215 | 27.52 | 112 | 23541 | 46.11 |
| 33 | 21345 | 19.32 | 73 | 21534 | 28.00 | 113 | 52341 | 46.89 |
| 34 | 25431 | 19.43 | 74 | 32154 | 29.65 | 114 | 41325 | 47.14 |
| 35 | 13245 | 19.55 | 75 | 43125 | 30.17 | 115 | 41532 | 47.50 |
| 36 | 35421 | 19.56 | 76 | 23451 | 30.40 | 116 | 53241 | 48.31 |
| 37 | 25314 | 19.97 | 77 | 24351 | 30.97 | 117 | 23415 | 48.84 |
| 38 | 45213 | 19.98 | 78 | 45123 | 31.07 | 118 | 42315 | 48.85 |
| 39 | 45312 | 20.07 | 79 | 24135 | 31.18 | 119 | 43215 | 52.91 |
| 40 | 21354 | 20.07 | 80 | 53412 | 31.22 | 120 | 41235 | 65.42 |

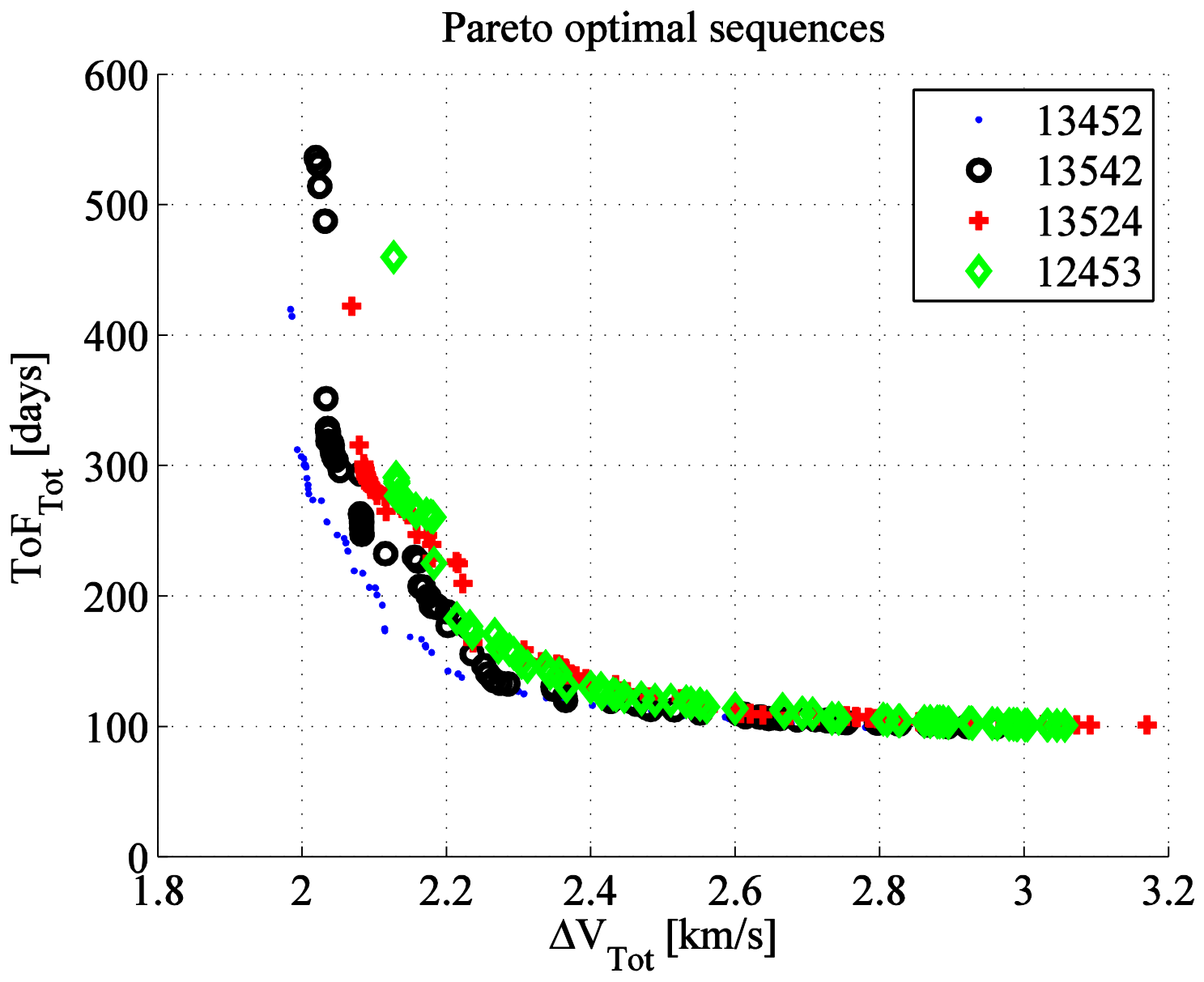

**Fig. 17** Pareto fronts corresponding to the four best sequences according to *Conv*

As one would expect, sequence 13452 has the lowest *Conv* since it coincides with part of the global Pareto front. Sequences 13524, 13542 and 12543 have also a low *Conv* index and thus they are quite close to the globally optimal solution, as shown in Fig. 17. In general, as already noted before, there is a strong dependence of the *quality* of the sequence from its first element. One can see that the first ranks are occupied mostly by sequences starting with debris nr. 1 and 3, while those with nr. 4 and 5 have highest *Conv* and are therefore occupy predominantly the worst ones. Those starting with nr. 2 are somewhat in the middle. The fact that solutions with nr. 1 and 3 are privileged as first elements in the sequence might be explained from the fact that they lie in the two lowest orbits (see Table 3) and therefore are *easier* to reach (Please keep in mind that for the rendezvous with the first debris there is no plane change since it is assumed to depart from a coplanar orbit). Another interesting observation is that the best sequences tend to avoid the largest plane changes. For example, in 13452 the plane changes are 1.47º, 2.52º, 1º and 2º. On the contrary, in the worst one according to *Conv*, i.e. 32415, they are 3.63º, 2.65º, 1.95º and 1 º.

**Table 7** Best $\Delta V_{Tot}$ and $ToF_{Tot}$ for each sequence**.** Best values are in **bold**. Worst values are underlined.

| **ord** | min($\Delta V_{Tot}$) [km/s] | min($ToF_{Tot}$) [days] | **ord** | min($\Delta V_{Tot}$) [km/s] | min($ToF_{Tot}$) [days] | **ord** | min($\Delta V_{Tot}$) [km/s] | min($ToF_{Tot}$) [days] |
|---|---|---|---|---|---|---|---|---|
| 12345 | 2.30 | 108.17 | 24513 | 2.17 | 106.64 | 42315 | 2.53 | 116.99 |
| 12354 | 2.26 | 107.81 | 24531 | 2.26 | 105.94 | 42351 | 2.63 | 114.24 |
| 12435 | 2.27 | 106.25 | 25134 | 2.24 | 102.96 | 42513 | 2.28 | 104.48 |
| 12453 | 2.13 | 100.63 | 25143 | 2.41 | 109.53 | 42531 | 2.36 | 109.50 |
| 12534 | 2.18 | 105.81 | 25314 | 2.30 | 107.66 | 43125 | 2.42 | 109.26 |
| 12543 | 2.13 | 103.26 | 25341 | 2.49 | 107.22 | 43152 | 2.34 | 107.36 |
| 13245 | 2.22 | 102.72 | 25413 | 2.26 | 104.36 | 43215 | 2.67 | 116.73 |
| 13254 | 2.11 | 103.03 | 25431 | 2.27 | 107.97 | 43251 | 2.56 | 113.93 |
| 13425 | 2.15 | 103.12 | 31245 | 2.20 | 100.73 | 43512 | 2.43 | 108.12 |
| 13452 | **1.98** | **96.35** | 31254 | 2.10 | 103.46 | 43521 | 2.53 | 111.63 |
| 13524 | 2.07 | 101.03 | 31425 | 2.30 | 106.79 | 45123 | 2.46 | 106.61 |
| 13542 | 2.02 | 100.08 | 31452 | 2.12 | 97.810 | 45132 | 2.33 | 102.12 |
| 14235 | 2.45 | 115.10 | 31524 | 2.15 | 104.32 | 45213 | 2.25 | 102.99 |
| 14253 | 2.21 | 104.45 | 31542 | 2.12 | 100.52 | 45231 | 2.42 | 110.85 |
| 14325 | 2.42 | 112.87 | 32145 | 2.44 | 111.30 | 45312 | 2.27 | 101.57 |
| 14352 | 2.21 | 105.30 | 32154 | 2.35 | 107.73 | 45321 | 2.55 | 111.39 |
| 14523 | 2.25 | 107.07 | 32415 | 2.51 | 115.31 | 51234 | 2.49 | 110.19 |
| 14532 | 2.27 | 105.43 | 32451 | 2.33 | 107.85 | 51243 | 2.38 | 107.21 |
| 15234 | 2.29 | 107.13 | 32514 | 2.38 | 107.70 | 51324 | 2.36 | 106.79 |
| 15243 | 2.14 | 102.56 | 32541 | 2.40 | 109.14 | 51342 | 2.25 | 103.38 |
| 15324 | 2.27 | 106.71 | 34125 | 2.42 | 109.68 | 51423 | 2.53 | 113.27 |
| 15342 | 2.17 | 102.45 | 34152 | 2.33 | 104.78 | 51432 | 2.55 | 112.77 |
| 15423 | 2.26 | 109.45 | 34215 | 2.36 | 112.19 | 52134 | 2.29 | 106.29 |
| 15432 | 2.24 | 106.63 | 34251 | 2.30 | 107.17 | 52143 | 2.44 | 108.36 |
| 21345 | 2.31 | 103.98 | 34512 | 2.18 | 101.86 | 52314 | 2.46 | 116.72 |
| 21354 | 2.24 | 103.07 | 34521 | 2.24 | 104.80 | 52341 | 2.61 | 112.97 |
| 21435 | 2.58 | 115.15 | 35124 | 2.27 | 103.81 | 52413 | 2.30 | 106.23 |
| 21453 | 2.38 | 106.26 | 35142 | 2.30 | 105.62 | 52431 | 2.37 | 108.24 |
| 21534 | 2.40 | 113.76 | 35214 | 2.32 | 109.19 | 53124 | 2.29 | 103.32 |
| 21543 | 2.32 | 110.97 | 35241 | 2.37 | 111.19 | 53142 | 2.36 | 108.06 |
| 23145 | 2.47 | 113.45 | 35412 | 2.28 | 101.50 | 53214 | 2.60 | 114.17 |
| 23154 | 2.36 | 107.94 | 35421 | 2.29 | 108.91 | 53241 | 2.62 | 116.59 |
| 23415 | 2.63 | 114.64 | 41235 | 2.70 | 116.91 | 53412 | 2.45 | 106.79 |
| 23451 | 2.48 | 111.27 | 41253 | 2.47 | 107.83 | 53421 | 2.46 | 112.91 |
| 23514 | 2.55 | 114.91 | 41325 | 2.55 | 113.15 | 54123 | 2.49 | 112.69 |
| 23541 | 2.54 | 111.22 | 41352 | 2.37 | 108.25 | 54132 | 2.38 | 105.83 |
| 24135 | 2.42 | 107.59 | 41523 | 2.54 | 111.57 | 54213 | 2.24 | 104.67 |
| 24153 | 2.43 | 108.50 | 41532 | 2.57 | 112.14 | 54231 | 2.42 | 115.17 |
| 24315 | 2.38 | 110.73 | 42135 | 2.44 | 115.59 | 54312 | 2.22 | 107.06 |
| 24351 | 2.42 | 108.28 | 42153 | 2.47 | 108.75 | 54321 | 2.50 | 116.42 |

Table 7 reports the minimum values for the performance parameters associated to each sequence, i.e. the extreme points of the Pareto fronts. Similar considerations to those made previously also apply to this case, with best values given by sequences starting with nr. 1 and 3 and the worst ones with nr. 4 and 5.

**Table 8** Debris removal sequence and timing for minimum $\Delta V_{Tot}$.

| Phase | *Final Keplerian elements* | | | | *Duration* | $\Delta V$ | *mass* [kg] |
|---|---|---|---|---|---|---|---|
| | *a* [km] | *e* | *i* [deg] | *Ω* [deg] | [days] | [km/s] | |
| Departure | 6628.16 | 0.010 | 1 | 65 | - | - | 1000 |
| Nr. 1 reached | 6828.16 | 0 | 1 | 65 | 5 | 0.115 | 996.11 |
| Nr. 1 de-orbited | 6752.69 | 0.011 | 1 | 65 | 22.06 | 0.043 | 993.21 |
| Nr. 3 reached | 6978.16 | 0 | -2 | 200 | 88.10 | 0.239 | 985.17 |
| Nr. 3 de-orbited | 6826.44 | 0.022 | -2 | 200 | 25.96 | 0.084 | 980.63 |
| Nr. 4 reached | 7478.16 | 0 | -1 | 90 | 66.71 | 0.476 | 964.88 |
| Nr. 4 de-orbited | 7055.54 | 0.053 | -1 | 90 | 34.33 | 0.221 | 951.69 |
| Nr. 5 reached | 7178.16 | 0 | 0 | 45 | 55.89 | 0.241 | 943.91 |
| Nr. 5 de-orbited | 6912.18 | 0.034 | 0 | 45 | 30.77 | 0.144 | 931.48 |
| Nr. 2 reached | 7128.16 | 0 | 2 | 150 | 56.98 | 0.297 | 922.12 |
| Nr. 2 de-orbited | 6901.39 | 0.032 | 2 | 150 | 33.99 | 0.124 | 917.24 |

Table 8 shows details about the best $\Delta V_{Tot}$ solution, with sequence 13452. Note that, in general, the $\Delta V$ cost of each phase is relatively low, thus leading to the minimum total cost of 1.98 km/s. Correspondingly, their duration is long, meaning that *slow* but more efficient transfers are preferred. This behaviour is also confirmed by the fact that the de-orbit conditions have non negligible eccentricities, which means also that the amplitude of the apogee thrusting arcs during de-orbit (see Fig. 2) is kept to a minimum. In this way propellant is devoted to lowering the perigee only with minimum variation of the apogee altitude.

By analysing in more detail the $\Delta V$ cost breakdown, one can see for example that the highest figures, 0.476 km/s are given by the rendezvous with debris nr. 4 from the de-orbit conditions of debris nr. 3. This high value is justified by the fact that reaching the final orbit radius of 7478.16 requires an apogee raise of 501 km from 6977 km and a perigee raise of 800 km from 6678 km. At the same time there is also a rotation of the orbit plane of 2.52°. By comparison, the rendezvous with nr. 5 after the de-orbit of nr. 4 is comparatively cheaper even if the radius of the target orbit is still high. In this case the perigee raise is 500 km while the apogee on the other hand needs to be lowered by 252 km from 7430 km since piece of debris nr. 4 is released on a relatively eccentric orbit with $e$=0.053. Plane rotation in this case is only 1°.

**Table 9** Debris removal sequence and timing for minimum $ToF_{Tot}$.

| Phase | *Final Keplerian elements* | | | | *Duration* | $\Delta V$ | *mass* [kg] |
|---|---|---|---|---|---|---|---|
| | *a* [km] | *e* | *i* [deg] | *Ω* [deg] | [days] | [km/s] | |
| Departure | 6628.16 | 0.010 | 1 | 65 | - | - | 1000 |
| Nr. 1 reached | 6828.16 | 0 | 1 | 65 | 5 | 0.115 | 996.11 |
| Nr. 1 de-orbited | 6685.24 | 0.001 | 1 | 65 | 4.04 | 0.081 | 990.61 |
| Nr. 3 reached | 6978.16 | 0 | -2 | 200 | 8.59 | 0.312 | 980.14 |
| Nr. 3 de-orbited | 6701.87 | 0.004 | -2 | 200 | 6.29 | 0.154 | 971.87 |
| Nr. 4 reached | 7478.16 | 0 | -1 | 90 | 14.79 | 0.664 | 950.17 |
| Nr. 4 de-orbited | 6789.06 | 0.016 | -1 | 90 | 17.13 | 0.362 | 928.76 |
| Nr. 5 reached | 7178.16 | 0 | 0 | 45 | 7.99 | 0.281 | 919.92 |
| Nr. 5 de-orbited | 6715.72 | 0.006 | 0 | 45 | 15.9 | 0.252 | 898.39 |
| Nr. 2 reached | 7128.16 | 0 | 2 | 150 | 9.87 | 0.466 | 884.28 |
| Nr. 2 de-orbited | 6725.90 | 0.007 | 2 | 150 | 6.75 | 0.221 | 875.87 |

Table 9 reports details about the minimum $ToF_{Tot}$ solution, again with sequence 13452. In contrast to what has been remarked for the previous case, here obviously the duration of each phase is kept to a minimum. For example, values for de-orbit times are very close to the minima reported in Table 3. Conversely, $\Delta V$ costs are higher than those in Table 8. Moreover, one can see that the de-orbit trajectories are quasi-circular, which suggests that the thrusting arcs are not restricted to apogee passages but cover almost entirely each revolution (i.e., with reference to Fig. 2, $\Delta L_a \approx 180°$).

A final note is devoted to the assumption mentioned in Section 3.2 that the delay due to phasing will be relatively negligible compared to the total transfer time. First of all, one has to consider that each de-orbit-rendezvous couplet is actually a transfer between two circular orbits with different altitude, phasing and orbit plane. In this sense, the related transfer strategy first lowers the perigee down to 300 km; then, in the second phase the apogee and perigee altitudes are adjusted to match those of the target orbit and at the same time the orbit plane is rotated around the line of nodes. In order to obtain a worst case estimation of the delay, it is chosen to decompose the latter into the contribution determined by the inclination change $t_{wait,\Delta i}$ and the one given by in-plane phasing $t_{wait,\Delta\phi}$. The former stems from the assumption made in Section 3.2 that the perigee lowering phase from the initial circular orbit is started such that the lines of apses is perpendicular to the line of nodes defined by the intersection of the orbit planes of the current piece of debris and the next one in the sequence. The maximum wait time is obtained when the line of nodes is aligned with the line of apses and is therefore given by half the orbit period of the departure circular orbit:

$$\max\left(t_{wait,\Delta i}\right) = \pi n_0 \quad (35)$$

where $n_0$ is the angular velocity of the initial circular orbit. After the line of apses is properly aligned in order to reach the target orbit plane, there remains, however, the problem of in-plane phasing. As a first step, the case of a quasi-circular transfer is considered, noting that this is actually the case for minimum $\Delta V$ sequences as the one reported in Table 8. If one considers the case of a transfer between two circular coplanar orbits, the phasing of the departure and arrival ones can by expressed as:

$$\begin{aligned} \phi_0 &= n_0 t_{wait,\Delta\phi} + \Delta\phi_{transf}\left(t_{transf}\right) \\ \phi_f &= n_f\left(t_{wait,\Delta\phi} + t_{transf}\right) + \Delta\phi_0 \end{aligned} \tag{35}$$

where $n$ is the angular velocity of the current orbit, $n_f$ is the one of the arrival orbit, is introducing along the transfer in order to match the phase of the arrival orbit. $\Delta\phi_0$ is the nominal phase difference between the two orbits at time of departure, computed simply from the initial and final argument of perigee and true anomaly:

$$\Delta\phi_0 = \left(\omega_f - \omega_0\right) + \left(\theta_f - \theta_0\right) \tag{35}$$

$\phi_0$ and $\phi_f$ can differ by multiples of $2\pi$, therefore, by combining Eqs. (35):

$$t_{wait,\Delta\phi}\left(n - n_f\right) = \Delta\phi_0 + \Delta\phi_{transf} + n_f t_{transf} + 2k\pi \quad k \in \mathbb{Z} \tag{35}$$

One can see that, once the transfer type is defined, the left side of Eq. (35) is constant and since $k$ is an arbitrary integer, one can write:

$$t_{wait,\Delta\phi} = \frac{\Delta\phi_{Tot}}{\left|n - n_f\right|} \quad \Delta\phi_{Tot} \in \left[0 \quad 2\pi\right] \tag{35}$$

and thus the worst case value for the delay $t_{wait}$ is obtained obviously for $\Delta\phi_{Tot} = 2\pi$. Since we are dealing with a LT transfer in which the semi-major axis is continuously varied, also the angular velocity $n$ at a certain point of the transfer is varying accordingly. Also, since it is assumed that the transfer is quasi-circular, one can insert a coasting arc of duration $t_{wait,\Delta\phi}$ at the point in which the ratio $1/\left|n - n_f\right|$ (which depends on the radii of the current and target orbits) is at its lowest. This condition typically occurs when the end of the de-orbit phase is reached.

If the transfer type is not quasi-circular but involves spirals with non negligible eccentricity, then an arbitrary delay cannot be introduced without altering the position of the lines of nodes. However, it is still possible to introduce an arbitrary number of coasting arcs of duration equal to the orbital period of the osculating orbit, i.e. one full revolution. The phase variation obtained by one such revolution is:

$$\Delta\phi_{2\pi}\left(n\right) = \frac{2\pi\left|n - n_f\right|}{n} \tag{35}$$

Note that, given the orbits involved in the transfer, $\Delta\phi_{2\pi}$ will be generally a fraction of $2\pi$. If a worst case phase variation $\Delta\phi_{Tot} = 2\pi$ is to be achieved, the following simple strategy can be used to estimate the corresponding delay: first, $k$ coasting revolutions are performed when the quantity $|n\text{-}n_f|$ is maximum. In this sense, one can write:

$$\begin{aligned} k &= \left\lfloor \frac{\Delta\phi_{2\pi}\left(n_k\right)}{2\pi} \right\rfloor \\ n_k &= \arg\min_n \Delta\phi_{2\pi}\left(n\right) \end{aligned} \tag{35}$$

This will bring the phase difference to a quantity which is lower than the maximum phase variation per revolution achievable, leaving a residual phase difference;

$$\Delta\phi_{res} = 2\pi - k\Delta\phi_{2\pi}\left(n_k\right) \tag{35}$$

A last coasting revolution is inserted to delete the residual when the semi-major axis which gives the proper angular velocity $n_{res}$ is reached:

$$n_{res} = \arg(\Delta\phi_{2\pi}\left(n\right) = \Delta\phi_{res}) \tag{35}$$

The total delay introduced in the worst case is therefore given by the sum of the periods of the coasting revolutions:

$$\max\left(t_{wait,\Delta\phi}\right) = 2\pi\left(\frac{k}{n_k} + \frac{1}{n_{res}}\right) \tag{35}$$

By applying the above strategies to the minimum $\Delta V$ and minimum time of flight sequences we can obtain a worst case estimation of the additional time introduced by phasing. The maximum delay introduced by the apses alignment in both cases would be 0.14 days. For the minimum time of flight case in Table 9 (i.e. quasi-circular sequence), the worst case delay due to $\Delta\phi$ is 2.68 days, leading to a total delay of 2.82 days. This value equates to a 2.93% increase compared to the nominal time of flight of 96.35 days, which can be considered acceptable for a preliminary study. On the contrary, in the case of minimum time of flight sequence as in Table 8, the delay due to $\Delta\phi$ would be 4.58 days, and the total delay 4.72 days, corresponding to a 1.12% increase on the nominal time of flight of 419.79 days. For these reasons, neglecting the phasing appears to be an acceptable approximation in this preliminary study.

## 5. Conclusions

This work presented a novel computational approach for the preliminary design of multispiral trajectories. The approach was applied to the design of an orbit debris removal mission by means of an IBS spacecraft. The models proposed here for the computation of low-thrust many-revolution transfers, allowed for a considerable reduction in control parameters and at the same time for a fast propagation of low-thrust motions thanks to the analytical propagation with FPET. Thanks to the reduced computational cost for the evaluation of a single fetch and deorbit operation, a multi-objective optimisation problem could be solved in which thousands of different debris removal sequences where examined to find the optimal ones with respect to $\Delta V$ cost and total removal duration. As a result, a considerable number of optimal candidate solutions where found. Analysis of the results showed that the particular removal sequence 13452 is globally optimal. A ranking criterion was proposed to grade all the candidate sequences and identify those that are suboptimal. From the analysis it was found that there is a dependency of the quality of the sequence on the first target object. Among the open issues for future developments, there is for example the possibility of integrating the problem of the sequence choice directly into a single multi-objective optimisation instance, thus obtaining a mixed continuous and discrete optimisation problem. This can be crucial when missions with more than 5-10 debris are considered since the decomposition in fixed-sequence continuous optimisation problems becomes less computationally tractable.

Future work will deal with comparing the proposed approach with similar methods like orbit averaging. In addition, although the proposed method has been applied to the special case of a debris removal mission, it is suitable to be extended and applied to more general trajectory design problems which involve many-revolution transfers from elliptical to circular or from elliptical to elliptical orbits Current developments are incorporating gravity perturbations in the analytical solution to allow the computation of more accurate solutions.

# 7. Appendix

In the following, the set of first-order solution for perturbed Keplerian motion. The Equinoctial elements at a longitude $L_f$ with respect to a reference longitude $L_0$ is given by a zero-order term plus a first order term multiplied by the magnitude $\varepsilon$ of the perturbing acceleration:

$$\begin{aligned}
a\left(L_f\right) &= a_0\left(L_0\right)+\varepsilon\, a_1\left(L_0,\Delta L,\alpha,\beta\right)\\
P_1\left(L_f\right) &= P_{10}\left(L_0\right)+\varepsilon\, P_{11}\left(L_0,\Delta L,\alpha,\beta\right)\\
P_2\left(L_f\right) &= P_{20}\left(L_0\right)+\varepsilon\, P_{21}\left(L_0,\Delta L,\alpha,\beta\right)\\
Q_1\left(L_f\right) &= Q_{10}\left(L_0\right)+\varepsilon\, Q_{11}\left(L_0,\Delta L,\alpha,\beta\right)\\
Q_2\left(L_f\right) &= Q_{20}\left(L_0\right)+\varepsilon\, Q_{21}\left(L_0,\Delta L,\alpha,\beta\right)\\
t\left(L_f\right) &= t_{00}\left(L_0,\Delta L\right)+\varepsilon\, t_{11}\left(L_0,\Delta L,\alpha,\beta\right)
\end{aligned} \tag{35}$$

where:

$$\Delta L = L_f - L_0 \tag{35}$$

For $a$, $P_1$, $P_2$, $Q_1$ and $Q_2$ the zero-order term is simply the value at $L_0$. For the time instead, is given by the reference time $t_0$ at $L_0$ plus the variation due to unperturbed Keplerian motion, where $I_{12}$ is an integral in $L$ as reported in Eq.(35).

$$\begin{aligned}
&a_0\left(L_0\right) = a_0\\
&P_{10}\left(L_0\right) = P_{10}\\
&P_{20}\left(L_0\right) = P_{20}\\
&Q_{10}\left(L_0\right) = Q_{10}\\
&Q_{20}\left(L_0\right) = Q_{20}\\
&t_{00}\left(L_0, L_f\right) = t_0 + h_0^{\,3} I_{12}
\end{aligned} \tag{35}$$

The first order terms are:

$$\begin{aligned}
a_1 &= 2\left(\frac{h_0 a_0}{\mu}\right)^2 \cos\beta\left(\cos\alpha\left(P_{20}I_{s2} - P_{10}I_{c2}\right) + \sin\alpha I_{11}\right)\\
P_{11} &= \frac{h_0^{\,4}}{\mu^3}\left(\cos\beta\left(-\cos\alpha I_{c2} + \sin\alpha\left(P_{10}I_{13} + I_{s2} + I_{s3}\right)\right) + \sin\beta P_{20}\left(-Q_{10}I_{c3} + Q_{20}I_{s3}\right)\right)\\
P_{21} &= \frac{h_0^{\,4}}{\mu^3}\left(\cos\beta\left(\cos\alpha I_{s2} + \sin\alpha\left(P_{20}I_{13} + I_{c2} + I_{c3}\right)\right) + \sin\beta P_{10}\left(Q_{10}I_{c3} - Q_{20}I_{s3}\right)\right)\\
Q_{11} &= 2\frac{h_0^{\,4}}{\mu^3}\sin\beta I_{s3}\left(1 + Q_{10}^{\,2} + Q_{20}^{\,2}\right)\\
Q_{21} &= 2\frac{h_0^{\,4}}{\mu^3}\sin\beta I_{c3}\left(1 + Q_{10}^{\,2} + Q_{20}^{\,2}\right)\\
t_1 &= \frac{\sqrt{a_0\left(1 - P_{10}^{\,2} - P_{20}^{\,2}\right)}}{2}\\
&\left(\begin{aligned}
&I_{c3}\left(2a_0\left(-3P_{10}P_{20}P_1 + P_2\left(2P_{10}^{\,2} - P_{20}^{\,2} - 2\right)\right) - 3a_1P_{20}\left(P_{10}^{\,2} + P_{20}^{\,2} - 1\right)\right)\\
&+I_{s3}\left(2a_0\left(-3P_{10}P_{20}P_2 + P_1\left(2P_{20}^{\,2} - P_{10}^{\,2} - 2\right)\right) - 3a_1P_{10}\left(P_{10}^{\,2} + P_{20}^{\,2} - 1\right)\right)\\
&-3I_{13}\left(2a_0\left(P_{10}P_1 + P_{20}P_2\right) + a_1\left(P_{10}^{\,2} + P_{20}^{\,2} - 1\right)\right)
\end{aligned}\right)
\end{aligned} \tag{35}$$

The terms $I_{11}$, $I_{12}$, $I_{13}$, $I_{c2}$, $I_{c3}$, $I_{s2}$, $I_{s3}$ are integrals in L, as:

$$I_{11}\left(L_f\right)=\int_{L_0}^{L_f}\frac{1}{1+P_{10}\sin L+P_{20}\cos L}dL$$

$$I_{12}\left(L_f\right)=\int_{L_0}^{L_f}\frac{1}{\left(1+P_{10}\sin L+P_{20}\cos L\right)^2}dL$$

$$I_{13}\left(L_f\right)=\int_{L_0}^{L_f}\frac{1}{\left(1+P_{10}\sin L+P_{20}\cos L\right)^3}dL$$

$$I_{c2}\left(L_f\right)=\int_{L_0}^{L_f}\frac{\cos L}{\left(1+P_{10}\sin L+P_{20}\cos L\right)^2}dL$$

$$I_{c3}\left(L_f\right)=\int_{L_0}^{L_f}\frac{\cos L}{\left(1+P_{10}\sin L+P_{20}\cos L\right)^3}dL$$

$$I_{s2}\left(L_f\right)=\int_{L_0}^{L_f}\frac{\sin L}{\left(1+P_{10}\sin L+P_{20}\cos L\right)^2}dL$$

$$I_{s3}\left(L_f\right)=\int_{L_0}^{L_f}\frac{\sin L}{\left(1+P_{10}\sin L+P_{20}\cos L\right)^3}dL \tag{35}$$

Finally Eqs. (35) to (35) report the complete analytical expressions for these integrals.

$$I_{11}=-\left.\frac{2\,\mathrm{atanh}\left(\dfrac{P10-\left(P20-1\right)\tan\left(\dfrac{L}{2}\right)}{\sqrt{P10^2+P20^2-1}}\right)}{\sqrt{P10^2+P20^2-1}}\right|_{L_0}^{L_f} \tag{35}$$

$$I_{12}=\left.-\frac{2\,\mathrm{atanh}\left(\dfrac{P10-\left(P20-1\right)\tan\left(\dfrac{L}{2}\right)}{\sqrt{P10^2+P20^2-1}}\right)}{\left(P10^2+P20^2-1\right)^{\frac{3}{2}}}+\frac{P10+\left(P10^2+P20^2\right)\sin\left(L\right)}{P20\left(P10^2+P20^2-1\right)\left(1+P20\cos\left(L\right)+P10\sin\left(L\right)\right)}\right|_{L_0}^{L_f} \tag{35}$$

$$I_{13}=\frac{1}{2}\left[-\frac{2\left(P10^2+P20^2+2\right)\operatorname{atanh}\left(\frac{P10-\left(P20-1\right)\tan\left(\frac{L}{2}\right)}{\sqrt{P10^2+P20^2-1}}\right)}{\left(P10^2+P20^2-1\right)^{\frac{5}{2}}}+\frac{P10+\left(P10^2+P20^2\right)\sin\left(L\right)}{P20\left(P10^2+P20^2-1\right)\left(1+P20\cos\left(L\right)+P10\sin\left(L\right)\right)^2}-\frac{P10\left(P10^2+P20^2+2\right)+3\left(P10^2+P20^2\right)\sin\left(L\right)}{P20\left(P10^2+P20^2-1\right)^2\left(1+P20\cos\left(L\right)+P10\sin\left(L\right)\right)}\right]\Bigg|_{L_0}^{L_f} \tag{35}$$

$$I_{c2}=-\frac{2P20\operatorname{atanh}\left(\frac{P10-\left(P20-1\right)\tan\left(\frac{L}{2}\right)}{\sqrt{P10^2+P20^2-1}}\right)}{\left(P10^2+P20^2-1\right)^{\frac{3}{2}}}-\frac{P10+\sin\left(L\right)}{\left(P10^2+P20^2-1\right)\left(1+P20\cos\left(L\right)+P10\sin\left(L\right)\right)}\Bigg|_{L_0}^{L_f} \tag{35}$$

$$I_{c3}=\frac{1}{2}\left[\frac{6P20\operatorname{atanh}\left(\frac{P10-\left(P20-1\right)\tan\left(\frac{L}{2}\right)}{\sqrt{P10^2+P20^2-1}}\right)}{\left(P10^2+P20^2-1\right)^{\frac{5}{2}}}-\frac{P10+\sin\left(L\right)}{\left(P10^2+P20^2-1\right)\left(1+P20\cos\left(L\right)+P10\sin\left(L\right)\right)^2}-\frac{3P10+\left(1+2P10^2+2P20^2\right)\sin\left(L\right)}{\left(P10^2+P20^2-1\right)^2\left(1+P20\cos\left(L\right)+P10\sin\left(L\right)\right)}\right]\Bigg|_{L_0}^{L_f} \tag{35}$$

$$I_{s2} = -\frac{2P10\,\mathrm{atanh}\left(\dfrac{P10-(P20-1)\tan\left(\dfrac{L}{2}\right)}{\sqrt{P10^2+P20^2-1}}\right)}{\left(P10^2+P20^2-1\right)^{\frac{3}{2}}} + \frac{P20^2-1+P10\sin(L)}{P20\left(P10^2+P20^2-1\right)\left(1+P20\cos(L)+P10\sin(L)\right)}\Bigg|_{L_0}^{L_f} \tag{35}$$

$$I_{s3} = \frac{1}{2}\left[\frac{6P10\,\mathrm{atanh}\left(\dfrac{P10-(P20-1)\tan\left(\dfrac{L}{2}\right)}{\sqrt{P10^2+P20^2-1}}\right)}{\left(P10^2+P20^2-1\right)^{\frac{5}{2}}} + \frac{P20^2-1-P10\sin(L)}{P20\left(P10^2+P20^2-1\right)\left(1+P20\cos(L)+P10\sin(L)\right)^2} \frac{P10\left(3P10+\left(1+2P10^2+2P20^2\right)\sin(L)\right)}{P20\left(P10^2+P20^2-1\right)^2\left(1+P20\cos(L)+P10\sin(L)\right)}\right]\Bigg|_{L_0}^{L_f} \tag{35}$$